\documentclass[11pt]{article}

\title{Electric Vehicle Scheduling with Capacitated Charging Stations and Partial Charging}
\date{\today}

\usepackage[natbibapa]{apacite}

\usepackage[a4paper, margin=1in]{geometry} 
\usepackage{graphicx,graphics,amsmath,amssymb, amsthm,
titling,graphicx,booktabs,tikz,subcaption,changepage,enumitem,rotating,lscape}

\usepackage{float}
\usepackage{graphicx}
\graphicspath{ {./Figures/} }
\usepackage[table]{}
\usepackage[ruled, lined, linesnumbered, commentsnumbered, longend]{algorithm2e}
\usepackage{ stmaryrd }
\usepackage{eurosym, outlines}
\usepackage{enumitem}
\usepackage{mathtools}
\usepackage{array}
\usepackage{booktabs}
\usepackage{rotating}
\usepackage{pifont}
\newcolumntype{P}[2]{%
  >{\begin{turn}{#1}\begin{minipage}{#2}\small\raggedright\hspace{0pt}}l%
  <{\end{minipage}\end{turn}}%
}

\usepackage{subcaption}

\usepackage{chngcntr}
\counterwithin{figure}{section}
\counterwithin{table}{section}
\usepackage{bbm}
\usepackage{setspace}
\usepackage[hidelinks]{hyperref} 

\newcolumntype{C}[1]{>{\centering\let\newline\\\arraybackslash\hspace{0pt}}m{#1}}
\newcolumntype{R}[1]{>{\raggedleft\let\newline\\\arraybackslash\hspace{0pt}}m{#1}}
\usepackage{dcolumn}
\usepackage[T1]{fontenc}
\usepackage{authblk}

\usepackage{authblk}

\begin{document}

    \author[1]{Marelot H. de Vos}
	\author[2,3]{Rolf N. van Lieshout}
	\author[4]{Twan Dollevoet}

	\affil[1]{\footnotesize ORTEC Data Science \& Consulting, Zoetermeer, The Netherlands}
	\affil[2]{\footnotesize Corresponding author. Email: r.n.v.lieshout@tue.nl}
	\affil[3]{\footnotesize Department of Operations, Planning, Accounting, and Control, School of Industrial Engineering, Eindhoven University of Technology, The Netherlands}
\affil[4]{\footnotesize Econometric Institute and ECOPT, Erasmus University Rotterdam, The Netherlands}

		\maketitle
	\begin{abstract}
    	\noindent This paper considers the scheduling of electric vehicles in a public transit system. Our main innovation is that we take into account that charging stations have limited capacity, while also considering partial charging. To solve the problem, we expand a connection-based network in order to track the state of charge of vehicles and model recharging actions. We then formulate the electric vehicle scheduling problem as a path-based binary program, whose linear relaxation we solve using column generation. We find integer feasible solutions using two heuristics: price-and-branch and truncated column generation, including acceleration strategies. We test the approach using data of the concession Gooi en Vechtstreek in the Netherlands, containing up to 816 trips. The truncated column generation outperforms the other heuristic, and solves the entire concession within 28 hours of computation time with an optimality gap less than 3.5 percent.

	\end{abstract}

\section{Introduction}\label{section:Introduction}
The benefits of electric buses are undisputed: replacing conventional combustion engine buses by electric buses drastically reduces noise, pollution and greenhouse gas emissions. For these reasons, many public transit operators have started electrifying their fleets. However, the introduction of electric buses introduces new complexities in the transit planning chain, since limited battery capacity requires electric buses to recharge during the day. In order to prevent buses from depleting their batteries while minimizing operating costs and energy consumption, it is crucial to time these recharging actions carefully, establishing a clear need for solution methods that incorporate charging to support bus companies in these decisions.

Naturally, the demand for algorithmic support for the planning of electric buses has spurred the interest in the Electric Vehicle Scheduling Problem (E-VSP), which is the problem of constructing feasible electric bus duties to cover a set of timetabled trips. A recent survey on electric bus planning and scheduling identified over 20 papers on (variants of) the E-VSP, primarily from the last three years \citep{perumal2021electric}. The innovation of our paper is that we propose a solution approach that is capable of solving large instances of the E-VSP while considering both the capacity of charging stations and partial (non-linear) charging, bridging the gap between theory and practice. 

A number of other papers on the E-VSP consider charging station capacity \citep{Li2014TransitEnergy,li2019mixed,tang2019robust,rinaldi2020mixed,wu2022multi}.  However, these papers either consider battery swapping or another form of constant-time charging. Another stream of literature considers partial charging \citep{Wen2016AnProblem,olsen2020scheduling,li2020joint,VanKootenNiekerk2017SchedulingVehicles,VanAken2020StrategischeOptimizer}. Literature on the E-VSP with both partial charging and capacitated charging stations is scarce \citep{Posthoorn2016VehicleApproach,janovec2019exact,zhang2021optimal}. \cite{janovec2019exact} only consider relatively small instances. \cite{Posthoorn2016VehicleApproach} tests the proposed approach on larger (single-depot) instances, but does not report optimality gaps or computation times. \cite{zhang2021optimal} solve medium-sized instances up to an optimality gap of 1\%, but consider only a single depot that also serves as the charging station. Conversely, in this paper, we find provably high-quality solutions for large instances with up to 816 trips of a rich variant of the E-VSP with partial charging, multiple capacitated charging stations, multiple depots and multiple vehicle types. 

Our solution approach is based on a discretization of the battery energy levels, which we combine with a connection-based network with nodes representing trips and charging actions. Every path in the resulting so-called \textit{primal} network represents a feasible bus duty, respecting both the compatibility of trips and the battery capacity. Note that the converse is not true: there may exist feasible vehicle duties that are not represented in the primal network, because we connect nodes using a conservative rounding scheme to ensure feasibility. However, we are able to use a relatively fine discretization, achieving a good trade-off between solution time and solution quality. In addition, although we do not consider this in our numerical experiments, the proposed discretization scheme perfectly lends itself to non-linear charging functions and is therefore widely applicable. 

Using the developed network structure, we formulate the problem as a path-based binary program with side constraints, whose linear relaxation can be solved with column generation. Due to the construction of the network, the pricing problem corresponds to a standard shortest path problem. To find integer solutions, we consider two heuristics: price-and-branch and truncated column generation. In price-and-branch, only the linear relaxation is solved using column generation, after which all generated paths are fed to a commercial MIP solver that optimizes over this given subset of all paths. Truncated column generation can be viewed as a diving heuristic in the branch-and-bound tree, where all nodes are solved using column generation: paths are fixed iteratively until the solution is integer feasible. We also develop acceleration strategies to reduce the computation time of the heuristics. 

In general, a disadvantage of optimization models based on discretization is that the obtained dual bound is only valid for the discretized representation of the problem, and is therefore not a true bound of the underlying problem. Inspired by \cite{Boland2017TheProblem}, we propose to find true lower bounds for the E-VSP by solving a linear relaxation on a \textit{dual} network, containing the same nodes as the primal network, but where nodes are connected using an \textit{optimistic} rounding scheme. Every feasible vehicle duty corresponds to a path in this modified network, such that this procedure yields a lower bound that is valid irrespective of the discretization. Note that where \cite{Boland2017TheProblem} apply this principle to time-discretized networks, we apply it to a network that is discretized in two dimensions: time and the battery energy levels. 

We test our approach using real-life timetable data from the bus concession Gooi en Vechtstreek in the Netherlands, which consists of 816 trips connecting five medium-sized cities southeast of Amsterdam. We also generate smaller instances by taking (random) subsets of all trips. On the smaller instances, truncated column generation achieves an optimality gap smaller than 1.5\%, outperforming price-and-branch. When paired with a dedicated acceleration strategy, we are able to solve the entire concession using truncated column generation up to an optimality gap of 3.4\%. We also perform a sensitivity analysis, which shows that our level of discretization is adequate: using finer discretizations leads to a significant increase in computation time, without resulting in a significant decrease in costs. 

The remainder of this paper is organized as follows: Section \ref{section:problem_description} gives a detailed description of the problem considered in this paper. Thereafter, Section \ref{section:literature_review} discusses literature related to this research. In Section \ref{section:metholodology}, we present our solution approach, discussing both the network structure and the column generation based heuristics. In Section \ref{section:results}, we present numerical results, including sensitivity analyses. Lastly, we conclude in Section~\ref{section:conclusion}.

\section{Problem Description}\label{section:problem_description}
The problem that we consider in this paper is a multi-depot vehicle scheduling problem that includes range constraints of the vehicles in a heterogeneous vehicle fleet, capacitated charging stations and partial charging. In this section, we discuss these elements in more detail. 

The core of the considered problem is the classical (multi-depot) vehicle scheduling problem (MD-)VSP \citep{Bunte2009AnModels}. Here, the input is a set of timetabled trips, which should all be performed by a single vehicle. Every trip has fixed starting and ending locations, as well as fixed starting and ending times. A pair of trips $(a,b)$ is called \textit{compatible} if their starting and ending times and locations are such that trip $b$ can be performed after trip $a$, potentially after an empty or \textit{deadhead} trip between the ending location of $a$ and the starting location of $b$. Every vehicle should be assigned a feasible duty, starting at (one of) the depot(s), performing a sequence of compatible trips, and returning to the (same) depot.

The problem considered in this paper is an extension of the MD-VSP with range constraints, (partial) recharging and capacitated charging stations. We assume that every electric vehicle is fully charged at the beginning of the day and that the energy consumption of every trip is known. Evidently, vehicles can only operate as long as their \textit{state of charge} (SoC) is strictly positive.\textit{ Charging actions} may be scheduled at specified charging stations with given capacities. We do not require that every charging action fully recharges a vehicle's battery, i.e. we allow for partial charging. The duration of a charging activity correlates positively with the increase in the SoC, according to a known \textit{charging function}. The capacity of charging stations implies that only a limited number of charging actions can be scheduled simultaneously at each charging station. In addition, we also consider a heterogeneous fleet, where the battery limit, charging function, and the energy consumption rate are allowed to differ per type of vehicle. Finally, we assume that the objective is to minimize the Total Costs of Ownership, including both fixed costs (investment costs for each vehicle) and operational costs (variable costs per kilometer to account for crew, energy consumption and maintenance).

\section{Literature Review}\label{section:literature_review}

For a general review on electric vehicle scheduling and related problems, we refer to \cite{perumal2021electric}. Here, we limit ourselves to discussing research on electric vehicle scheduling that considers the capacity of charging stations and/or partial charging. An overview of the included features in the discussed literature is presented in Table~\ref{tab:overviewE-VSP}.

\begin{table}[ht]
  \centering
  \caption{Overview of included aspects in papers on the E-VSP. Abbrevations: CG=Column Generation (incl. branch-and-price and CG-based heuristics), MIP=Mixed Integer Programming, MH=Metaheuristic.}
  \label{tab:overviewE-VSP}
  \scalebox{0.80}{
  \begin{tabular}{r|cccccD{.}{.}{3.2}}
  \hline
  & 
     \multicolumn{1}{P{60}{1.6cm}}{Partial Charging} &
    \multicolumn{1}{P{60}{1.5cm}}{Multiple Depots} &
    \multicolumn{1}{P{60}{1.5cm}}{Multiple Vehicle Types} &
    \multicolumn{1}{P{60}{1.6cm}}{Charging Station Capacity} &
        \multicolumn{1}{P{60}{1.5cm}}{Solution Method} &
    \multicolumn{1}{P{60}{1.5cm}}{Number of Trips} \\
  \midrule
    \citet{Li2014TransitEnergy} 						&   &  &  &  \ding{51} & CG &  947\\
    \citet{Posthoorn2016VehicleApproach}  				&	\ding{51}   &  & & \ding{51}  & CG & 709 \\
    \citet{Wen2016AnProblem} 							& \ding{51} & \ding{51} &  &  & MH & 500\\
    \citet{VanKootenNiekerk2017SchedulingVehicles}  	& \ding{51} &   &  &  & CG  & 543\\
    \citet{janovec2019exact} 							& \ding{51} &  &  & \ding{51} & MIP  & 160\\
    \citet{li2019mixed} 								&   & \ding{51} & \ding{51}   & \ding{51} & MIP & 288\\
    \citet{tang2019robust} 								&   &  &  & \ding{51} & CG & 96\\
    \citet{li2020joint} 								& \ding{51} & \ding{51} &  & & MH  & 867\\
    \citet{olsen2020scheduling} 						& \ding{51} & \ding{51} &  & & MH  & 10{\small,}710\\
    \citet{rinaldi2020mixed} 							&   &  & \ding{51} & \ding{51} & MIP & 1008\\
    \citet{VanAken2020StrategischeOptimizer} 			& \ding{51}   & \ding{51}  & \ding{51} & & CG  & 1{\small,}200\\
    \citet{zhang2021optimal} 			& \ding{51}  &   &  & \ding{51} & CG  & 160\\

     \citet{wu2022multi} 			&   & \ding{51}  &  & \ding{51} & CG  & 400\\
      This paper 										& \ding{51}   & \ding{51}  & \ding{51} & \ding{51} & CG  & 816\\
  \hline
  \end{tabular}}
\end{table}

\subsection{Capacitated Charging Stations}
 \citet{Li2014TransitEnergy}  studies the E-VSP with battery swapping (or equivalently, fast-charging), taking the capacity of the charging station into account. In this setting, charging a vehicle takes a constant time. To solve the problem, the author develops exact and heuristic algorithms based on column generation. \cite{li2019mixed} deals with scheduling a mixed fleet of electric and conventional buses. It is assumed that full energy is restored after each refueling, which takes a constant time of 30 minutes, i.e. partial charging is not considered. The authors formulate the problem as a mixed-integer programming (MIP) model and solve instances with 288 trips and two depots using a commercial solver. \cite{tang2019robust} investigate the scheduling of electric buses in a stochastic setting, where the aim is to find schedules robust against varying traffic conditions. The authors include charging station capacity, but only consider fast charging. Using branch-and-price, problem instances with up to 96 trips are solved, both in static and in dynamic fashion. \cite{rinaldi2020mixed} propose a MIP model to schedule a mixed-fleet of electric and diesel buses, assuming fast-charging. The authors also develop an ad-hoc decomposition scheme, which is tested on instances with up to 1008 trips. \cite{wu2022multi} propose a branch-and-price scheme to solve the E-VSP with capacitated charging stations, minimizing both costs and the overall peak load on the energy grid. The authors use the epsilon-constraint method to find (approximate) Pareto-efficient solutions with respect to the two objectives for instances with up to 400 trips. 
 
\subsection{Partial Charging}
\cite{Wen2016AnProblem} develop a MIP model and an Adaptive Large Neighborhood Search heuristic to solve the E-VSP with partial charging. The MIP can solve instances with 30 trips, and the heuristic with up to 500 trips. \cite{olsen2020scheduling} extend this heuristic to allow for non-linear charging processes and analyze the impact of assuming a constant charging time and/or a linear charging process on large instances with thousands of trips. An alternative heuristic approach is taken by \cite{li2020joint}, who develop an adaptive genetic algorithm, which is used to solve instances with up to 867 trips. \cite{VanKootenNiekerk2017SchedulingVehicles} present two MIP models for the E-VSP with partial charging. The first model assumes a linear charging process, so that the SoC can be tracked with continuous variables. In the second model this assumption is relaxed, which requires the SoC to be discretized. The authors apply column generation based heuristics to solve instances with 543 trips using the second model. The first model is only solvable for small instances. A similar discretization approach is taken by \cite{VanAken2020StrategischeOptimizer}, who, besides partial charging, also consider multiple depots and bus types, and are able to solve instances containing up to 1,200 trips.

\subsection{Capacitated Charging Stations and Partial Charging}

\cite{Posthoorn2016VehicleApproach} considers the E-VSP with partial charging and a single charging station with a limited capacity. The author discretizes the SoC and solves the linear relaxation using column generation. Subsequently, the generated paths are included in a MIP to find integer solutions. The approach is applied to instances with up to 709 trips. The discretization is relatively coarse and no gaps or computation times are reported. 
 \cite{janovec2019exact} develop an arc-based MIP formulation for the E-VSP with partial charging and capacitated charging stations. The authors consider instances with up to 160 trips with nine buses and three to six chargers. Furthermore, for all instances, the results with electric buses are the same as with diesel buses, which suggests that the range and charging capacity constraints may not be restrictive. \cite{zhang2021optimal} study electric vehicle scheduling with partial (non-linear) charging from a single terminal, which also serves as the (capacitated) charging station. In addition, the authors also consider the impact of the schedule on battery-aging. Instances with up to 160 trips are solved using a branch-and-price algorithm.

\section{Methodology}\label{section:metholodology}
In order to solve the E-VSP, we first formulate the problem as a set covering problem with additional constraints. We then develop a column generation algorithm to solve its LP-relaxation. Finally, we present two heuristics, based on column generation, to obtain feasible solutions. 

\subsection{Mathematical Formulation}\label{section:meth-model} 

Our mathematical formulation is based on the set-partitioning model for the VSP as explained by \citet{Bunte2009AnModels}. In this formulation, the columns correspond to feasible vehicle duties, also called \textit{paths}, which are sequences of compatible trips. Given that we consider electric vehicles, these paths include recharging actions as well. In our formulation, we consider set covering instead of set partitioning, since the set covering formulation is claimed to be numerically more stable \citep{Barnhart1998IntegerPrograms}. Hence, each trip should be serviced by at least one, instead of precisely one vehicle. Without loss of generality, a double trip can be deleted from a vehicle duty without increasing the costs, and therefore an optimal solution to the set covering formulation is also optimal for the set partitioning formulation. 

Moreover, since we solve the E-VSP, charging activities should be taken into account. Similar to \citet{Li2014TransitEnergy}, we discretize the time horizon into time blocks $\mathcal{B}$ with a fixed length $l$. These time blocks are created to track the availability of charging stations over time and to incorporate the limited capacity of the charging stations. We assume that a vehicle occupies the charging station either during the entire block, or not at all in this block. Let $t_b$ represent the starting time of time block $b \in \mathcal{B}$. A time block $b \in \mathcal{B}$ represents the time interval $[t_b,t_b+l)$.

We can now formulate the E-VSP. 
The set $\mathcal{T}$ represents all trips that should be serviced, while the set $\mathcal{R}$ contains all charging stations. The parameter $M_r$ denotes the capacity of charging station $r \in \mathcal{R}$. We define the set $\mathcal{P}$ containing all possible paths. A path $p \in \mathcal{P}$ encodes a feasible vehicle duty in which trips to service and charging actions are included. We define a binary decision variable $x_p$ that indicates whether path $p \in \mathcal{P}$ is selected in the solution. The parameter $c_p$ represents the costs of path $p \in \mathcal{P}$, and each coefficient $a_{i,p}$ is 1 if trip $i \in \mathcal{T}$ is included in path $p \in \mathcal{P}$, and 0 otherwise. Similarly, the coefficient $u_{r,b,p}$ is 1 if charging station $r \in \mathcal{R}$ is visited during time block $b \in \mathcal{B}$ in path $p \in \mathcal{P}$. We use the following formulation for the E-VSP:
\begin{align}
    \min \quad & \sum_{p\in \mathcal{P}} c_p x_p, \label{eq:setPart-obj}\\ 
   \textnormal{s.t.}  \quad & \sum_{p\in \mathcal{P}} a_{i,p} x_p \geq 1 &\forall\ i \in \mathcal{T}, \label{eq:setPart-trip}\\
   & \sum_{p\in \mathcal{P}} u_{r,b,p} x_p \le M_r &\forall\ r \in \mathcal{R}, \ b \in \mathcal{B} ,  \label{eq:setPart-capacity}\\
    & x_{p} \in \{0,1\} & \forall p \in \mathcal{P}. \label{eq:setPart-binary}
\end{align}
The Objective \eqref{eq:setPart-obj} minimizes the total costs. Constraint \eqref{eq:setPart-trip} guarantees that each trip $i \in \mathcal{T}$ is executed by at least one vehicle. Constraint \eqref{eq:setPart-capacity} is added to ensure that the capacity of each charging station $r \in \mathcal{R}$ is not exceeded in any of the time blocks $b \in \mathcal{B}$. Lastly, Constraint \eqref{eq:setPart-binary} provides the range of the decision variables.

\subsection{Column Generation Algorithm}

Since the path-based formulation \eqref{eq:setPart-obj}-\eqref{eq:setPart-binary} has exponentially many variables, we develop two heuristics that are based on column generation. In column generation, a solution to a linear program, called the \textit{master problem} (MP), is found by iteratively solving a \textit{restricted master problem} (RMP) and a \textit{pricing problem}. In short, the MP is the LP-relaxation of the set covering model \eqref{eq:setPart-obj}-\eqref{eq:setPart-binary}. The RMP is similar to the MP but uses only a subset of paths, denoted by $\mathcal{P'}$. Before the start of the column generation process, the set $\mathcal{P'}$ is initialized with paths that allow a feasible solution to the RMP. Afterward, in every iteration, the RMP is solved and the values of the dual variables are used as input for the pricing problem. The pricing problem searches for new paths with negative reduced costs, since these can improve the objective value. The path with the most negative reduced costs is added to the set $\mathcal{P'}$ used in the RMP. If the pricing problem cannot provide a path with negative reduced costs, the MP is solved to optimality and the column generation process is terminated. For a more detailed explanation of column generation, we refer to \citet{Desrosiers2005AGeneration}. 

In the remainder of this section, we discuss the pricing problem of the column generation approach for the E-VSP in detail. Additionally, we discuss how the set $\mathcal{P}'$ is initialized. 
We then discuss how lower bounds on the optimal solution value can be obtained.
Finally, we develop two column generation based heuristics that are used to obtain a feasible solution for the E-VSP.

\subsubsection{Pricing Problem}\label{section:meth-PricingProblem}

In every iteration, the pricing problem searches for new variables with negative reduced costs based on the current values of the duals. From the RMP, we obtain optimal values for the dual variables $\sigma_i$ for all trips $i \in \mathcal{T}$ and $\gamma_{r,b}$ for all charging stations $r \in \mathcal{R}$ and time blocks $b \in \mathcal{B}$. Using this dual information from the RMP, the reduced costs corresponding to path $p$ can be calculated as follows:
\begin{align}
    RC(x_p) = c_p - \sum_{i \in \mathcal{T}}a_{i,p}\sigma_i - \sum_{r \in \mathcal{R}}\sum_{b \in \mathcal{B}} u_{r,b,p}\gamma_{r,b}.
\end{align}

The aim of the pricing problem is to find the variable corresponding to a path with the lowest reduced costs. In order to find the path $p\in \mathcal{P}$ with the lowest reduced cost, we solve a shortest path problem in a suitably chosen network, which we now describe in full detail.

\paragraph{Network Structure}\label{section:meth-network}
We construct a connection-based network which allows to create feasible vehicle duties. In particular, we create a separate network for each combination of a vehicle type and depot, as proposed by \citet{Gintner2005SolvingPractice}. This allows us to ensure that each vehicle starts and ends at the same depot and to incorporate different characteristics for each vehicle type. We define $\mathcal{K}$ as the set of networks that are created. The aim is to find a vehicle duty with lowest reduced costs in each network $k \in \mathcal{K}$ separately.

We now describe the construction of the network $G^{k}(\mathcal{N}^{k},\mathcal{A}^{k})$ for a given combination $k \in \mathcal{K}$ of a vehicle type and a depot. We take the SoC of the vehicle into account by discretizing the possible SoC values and tracking the SoC of the vehicle along the path. Thus, the nodes in this network represent the depot, combinations of trips and SoC values, or combinations of charging stations, time blocks, and SoC values. Each compatible connection is explicitly modeled using a conservative rounding scheme for the SoC values. This ensures that all paths in the network correspond to a feasible vehicle duty. Therefore, we refer to this network as the \textit{primal} network. However, as a consequence, some feasible vehicle duties cannot be represented as a path in our network. we introduce the concept of a \textit{dual} network and this network can generate true lower bounds in Section~\ref{section:lowerBounds}. Furthermore, we study the impact of the discretization on the solution values in Section~\ref{section:results-discretization}. 

\paragraph{Nodes}
In each network $G^k$, we include a source node and sink node, denoted as $d^{\sigma}_k$ and $d^{\tau}_k$, respectively, that correspond to the depot. We include nodes for the timetabled trips in combination with discretized SoC values similar to \citet{VanKootenNiekerk2017SchedulingVehicles}, to keep track of the SoC of vehicles along their path. Mathematically, let $\mathcal{T}^{k}$ and $\mathcal{S}^{k}$ be the sets of trips and all possible SoC values for network $k$, respectively. Since we discretize the SoC, $\mathcal{S}^{k}$ has a finite number of elements. Let $s^\textnormal{min}_{k}$ represent the minimum allowed SoC value. For trip $i$, we define the set $\mathcal{S}^{k}_{i}$ that includes all SoC values $s \in \mathcal{S}^{k}$ that are at least $s^\textnormal{min}_{k}$ plus the SoC required for executing trip $i$, represented by $f_i$. This results in the node set $$\mathcal{N}_{k}^\textnormal{trip}=\{(i,s)|i \in \mathcal{T}^{k}, s \in \mathcal{S}^{k}_{i}\}.$$ 
For each node, the value $s$ represents the SoC of the vehicle at the moment it departs from the start location of trip $i$, hence at the beginning of trip $i$.

Additionally, to be able to model charging activities, we use copies of charging actions in combination with possible SoC values as nodes. We mean by \textit{charging action} the charging of a vehicle at a charging station within a specific time block $b\in \mathcal{B}$. Let the set $\mathcal{R}^k$ correspond to all charging stations suitable for network $k$.
As the SoC value at a specific node influences the compatibility to other nodes in the network, we combine these nodes with discretized SoC values. We include nodes per charging action for each possible value of the SoC $s \in \mathcal{S}^{k}$ that a vehicle has at the beginning of the charging action for all charging nodes. Since this SoC value represents the SoC before the charging, we disregard the nodes for a fully charged SoC, represented by $s^\textnormal{full}$, since this value cannot be increased during a charging action. Thus, the nodes created for the charging actions can be summarized in the set 
$$\mathcal{N}_{k}^\textnormal{charge}=\{(r,b,s)|r \in \mathcal{R}^k, b \in \mathcal{B}, s \in \mathcal{S}^{k}\setminus \{s^\textnormal{full}\}\}.$$ 

In conclusion, the nodes in network $k$ are represented by the node set $$\mathcal{N}^{k}=\{d^{\sigma}_{k}, d^{\tau}_{k}\} \cup \mathcal{N}_{k}^\textnormal{trip} \cup \mathcal{N}_{k}^\textnormal{charge}.$$ 
For the remainder of this paper, we call node $n$ a \textit{trip node} if $n \in \mathcal{N}_{k}^\textnormal{trip}$ and a \textit{charging node} if $n \in \mathcal{N}_{k}^\textnormal{charge}$ in a particular network $k \in \mathcal{K}$. 

\paragraph{Arcs}

The set of arcs $\mathcal{A}^k$ represent connections between nodes in the network. Because the time horizon and SoC values are discretized, a rounding scheme is needed. For bus companies, it is important that an electric bus is always capable of finishing its duty and, thus, running out of battery must be prevented. Additionally, a bus arriving too early is preferred over that bus arriving too late. For these reasons, we apply a conservative rounding scheme. In particular, charging actions begin at the earliest in the first time interval $b \in \mathcal{B}$ that starts after the arrival of the vehicle. The vehicle idles in between its arrival and the start of the time block. Additionally, we round down the actual SoC value of a vehicle to the nearest SoC value $s \in \mathcal{S}^{k}$, which results in an underestimated SoC value. Consequently, it can occur that some feasible vehicle duties are excluded. However, the final schedule obtained using this conservative rounding scheme will certainly be feasible in practice.

Below, we briefly explain which arcs are created in network $k$. Here, we let $F^\textnormal{soc}(\cdot)$ be a function that returns the nearest SoC value $s \in \mathcal{S}^{k}$ smaller than the SoC input. A more thorough description is provided in Appendix~\ref{appendix:arcs}.
\begin{itemize}
	\item Given that vehicles leave the depot fully charged, we only include arcs from the source node $d^\sigma_k$ to trip nodes. For these arcs, we take the SoC usage required for the deadheading trip from the depot to the start location of the corresponding trip into account, as well as a maximum deadheading time, if applicable.
	\item For trip nodes $n=(i,s)$, we include outgoing arcs to the sink node, to other trip nodes, and to charging nodes. An outgoing arc to the sink node is included if the node's SoC value is sufficient to execute the trip and then directly return to the depot. An outgoing arc to another trip node $v$ is present if the trips are compatible and the SoC value $s$ is sufficient to execute both trips, as well as the deadheading and idling between the trips, if applicable. Similarly, there is an outgoing arc from a trip node to a charging action $v$ if the maximum deadheading and idling SoC usage and time are respected. For arcs towards a trip or charging node $v$, the SoC value $s'$ of $v$ must satisfy
	$$s'= F^\textnormal{soc}(s - f_{i} - \tau^\textnormal{soc}_{(n,v)}),$$
	where $\tau^\textnormal{soc}_{(n,v)}$ represents the SoC required for the potential deadheading and idling between nodes $n$ and $v$, if applicable.
	\item For charging nodes $n=(r,b,s)$, we include arcs to the sink node, to trip nodes, and to other charging nodes. We denote the amount of SoC recharged in node $n$ by $s_n^+$. An outgoing arc to the sink node is included if the node's SoC value after recharging is sufficient to return to the depot, while this is not the case without recharging $s_n^+$. If it would be possible to return to the depot without the recharging of node $n$, vehicles could either return to the depot without recharging at all, or with a shorter recharging time. Outgoing arcs to trip nodes $v$ are present if the vehicle can arrive timely at the start location of the trip. Moreover, the SoC value $s'$ of the trip node must equal $s'= F^\textnormal{soc}(s + s_n^+ - \tau^\textnormal{soc}_{(n,v)})$. Finally, an outgoing arc from one charging node to another is included if the corresponding charging stations coincide, the second charging node directly follows the first one in time, and $s'= F^\textnormal{soc}(s + s_n^+) > s$.
\end{itemize}

\paragraph{Arc Costs}
The costs of a path $p\in \mathcal{P}$ can be distributed over the arcs. Each arc $(n,v) \in \mathcal{A}^{k}$ contains operational costs per driven kilometer and crew costs for deadheading and idling. These general costs for arcs are represented by $c^\textnormal{arc}_{(n,v)}$. Arcs coming from the source node ($n = d^{\sigma}_{k}$) additionally include the investment costs corresponding to the vehicle of network $k$, represented by $c^\textnormal{invest}_{k}$. Each arc that reaches a trip node includes operational costs per driven kilometer and the crew costs of the corresponding trip. Similarly, each arc that reaches a charging node, includes the operational costs of the corresponding charging action.
Each arc that reaches a charging node coming from a trip node, includes an additional fixed penalty term $c^\textnormal{\textnormal{start}}$, which can be interpreted as starting costs for a charging activity. This penalty term should be small relative to the investment costs. 

The reduced costs of a path $p\in\mathcal{P}$ can be distributed over the arcs as well. Besides the primal costs that are described above, we subtract the dual costs $\sigma_i$ for all trip nodes $v = (i,s)\in\mathcal{N}_{k}^\textnormal{trip}$ from the arc costs of all arcs $(n,v)$ towards $v$, and we subtract the dual costs $\gamma_{r,b}$ for all charging nodes $v = (r,b,s) \in\mathcal{N}_{k}^\textnormal{charge}$ from the arc costs of all arcs $(n,v)$.

\begin{figure}[ht]
    \centering
    \includegraphics[scale=0.7]{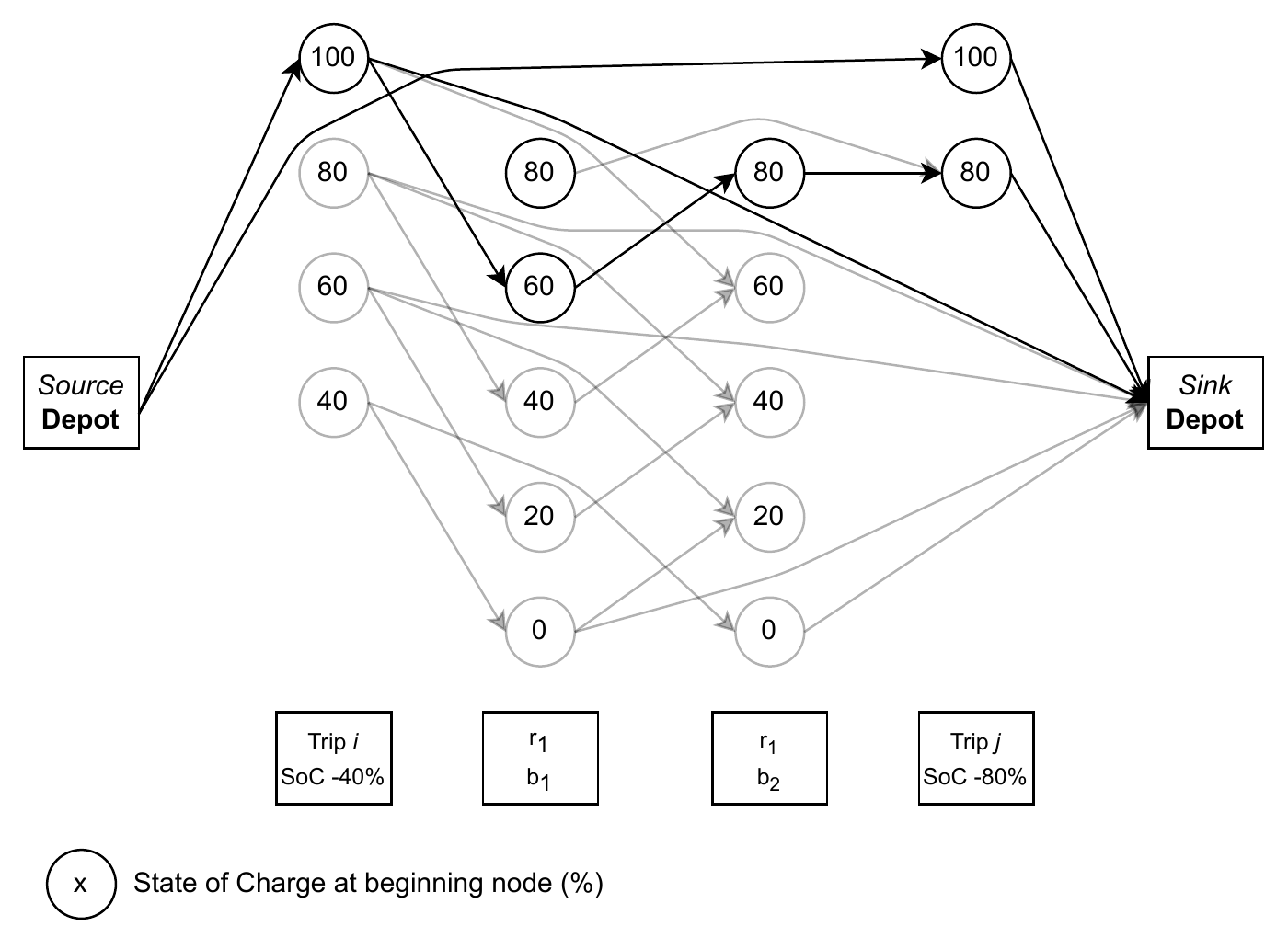}
    \caption{Network including all nodes and arcs before the preprocessing. The nodes and arcs with reduced opacity are not connected to the source node and/or sink node, and are removed in a preprocessing step.}
    \label{fig:voorbeeldnetwork}
\end{figure}

\paragraph{Example} An example of a resulting network is given in Figure~\ref{fig:voorbeeldnetwork}. Here, we have two trips $i$ and $j$ that can be serviced by the vehicle type under consideration and a charging station $r_1$ that is available during two time blocks $b_1$ and $b_2$ between the end of trip $i$ and the start of trip $j$. Both trips start at the depot. Trip $i$ and $j$ reduce the SoC by 40\% and 80\%, respectively. Trip $i$ ends at a location from which the SoC is reduced by 20\% to return to the depot. Trip $j$ ends at the depot. Charging station $r_1$ is located at the end location of trip $i$ and can increase the SoC in one time block with 20\%, independently of the SoC value at the beginning of the charging action. 
We use SoC values between 0\% to 100\% in steps of 20\%. 
Observe that some of the nodes and arcs cannot be included in a path that starts and ends at the depot. These nodes and arcs are depicted in a lighter color. In a preprocessing step, we remove these node and arcs from the network.

\subsubsection{Initialization}
To initialize the column generation algorithm, we construct a feasible solution by including a separate path for each trip in $\mathcal{P}'$. In our test instances, this solution is always feasible and no charging is necessary. In general, one can also avoid having to find a set of feasible paths for the initialization by using dummy variables in the trip covering constraints. 

\subsubsection{Lower bounds}\label{section:lowerBounds}
 The column generation algorithm terminates if paths with negative reduced costs can no longer be found. In that case, the MP has been solved to optimality, given the set $\mathcal{P}$ of columns under consideration. A lower bound on the MP's solution value can be obtained in each iteration of the column generation algorithm, if a value $\kappa$ can be found that satisfies
$$\kappa \geq \sum_{p \in \mathcal{P}} x_p$$
for all optimal solutions of the MP. In that case, the value $z_{\textnormal{RMP}} + \kappa c$ is a lower bound on the optimal solution value of the MP, where $c$ is the lowest reduced costs over all paths $p \in \mathcal{P}$ and $z_{\textnormal{RMP}}$ is the optimal solution value of the RMP (see \cite{Desrosiers2005AGeneration}).

As stated before, because of the conservative rounding scheme, there might be feasible vehicle duties that cannot be represented as a path in the networks defined above. Hence, the lower bound on the MP does not necessarily correspond to a true lower bound of the original problem. 

To obtain a true lower bound, we construct an auxiliary dual network that has the same nodes as the primal network, but where nodes are connected using an optimistic rounding scheme instead. Such a scheme rounds SoC values up and allows to start charging in the last time block before the vehicle arrives at a charging station. It should be noted that not all paths in the dual network correspond to feasible vehicle duties. Any lower bound on the MP's solution value formulated on the dual network corresponds to a lower bound that is valid irrespective of the used discretization. We compute such a lower bound as a separate step from finding feasible solutions.

\subsubsection{Obtaining Integral Solutions}\label{section:meth-heuristics}
We now discuss two heuristics to find integral solutions for the E-VSP. In both cases, we first obtain a solution to the master problem by using column generation.

\paragraph{(Truncated) Price-and-Branch} Our first heuristic solves the MP, and then applies a commercial solver to the binary program \eqref{eq:setPart-obj}-\eqref{eq:setPart-binary} using all columns in $\mathcal{P}'$. Such a heuristic has been referred to as a \textit{restricted master heuristic} or \textit{price-and-branch} \citep{Sadykov2019PrimalMethods}. For this heuristic, we can either solve the MP to optimality, or terminate the column generation algorithm before an optimal solution to the MP has been obtained. In the latter case, we prevent the infamous tail-off effect of column generation. We terminate the column generation process if the objective value of the RMP has not improved sufficiently relative to the objective value of a fixed number of iterations earlier. Let the parameters $Z_\textnormal{min}$ and $I$ represent the minimum percentage of improvement considered as sufficient and the number of iterations, respectively. Several existing papers use this method to early stop a column generation process \citep[see e.g.][]{Wang2019ColumnProblem,Gamache1999ColumnProblems}. After the column generation process is terminated, we solve the binary program as in the price-and-branch heuristic. If the column generation process is terminated early, we refer to this heuristic as \textit{truncated price-and-branch}.

\paragraph{Truncated Column Generation}
Our second heuristic is comparable to the truncated column generation as proposed by \citet{Pepin2009AProblem} to solve the MD-VSP. In this heuristic, a column generation phase and a fixing phase are executed iteratively. Specific path variables in the fractional solution are fixed to one after early termination of the column generation phase. Afterward, the column-generation phase re-starts. This iterative process is repeated until the obtained solution is integral. This is in contrast to the first proposed heuristic that applies a commercial solver after the column generation process to obtain an integer solution. We use the same stopping criterion for each column generation process as for the truncated price-and-branch using the predefined parameters $Z_\textnormal{min}$ and $I$. This early termination differs from \citet{Pepin2009AProblem} since we look at a relative difference, whereas \citet{Pepin2009AProblem} use an absolute difference. We use a relative difference to be able to apply a similar method to instances of varying sizes. After termination of each column generation phase, all path variables $x_p$ in the set $\mathcal{P}'$ that have a solution value larger than a predefined threshold $\theta$ and that are no initial path are set to one in the RMP if the solution is fractional. If no such paths exist, we fix the path with the maximum value that is not an initial path. Then, we resolve the RMP and start the column generation process again. Each time the column generation process is restarted, at least $I$ iterations should be executed before we early terminate the process. These steps are repeated until an integral solution is obtained. 

To reduce the network size during the column generation process, we can extend the truncated column generation heuristic by executing an additional operation each time paths are fixed. This additional step aims to lower the computation time for solving the pricing problem. In this operation, certain nodes and their adjacent arcs are removed from the pricing networks. Firstly, the nodes corresponding to all trips included in the paths that are fixed are removed. Choosing new paths that also contain these trips results in executing trips more than once. It is unlikely that this gives an optimal solution. Hence, it is reasonable to delete these nodes and their adjacent arcs.
Secondly, the nodes corresponding to each combination of a charging station and time block that reached the capacity of the charging station due to fixed paths can be deleted. Since the capacity constraints of these charging station and time block combinations are binding, newly added paths that include a charging action for such a combination can never be part of the solution. 
Lastly, the preprocessing step is repeated. In this step, all nodes that are unreachable from the source or cannot reach the sink are removed.
If all nodes as described above and their adjacent arcs are deleted, the network size shrinks once paths are fixed.  A potential disadvantage of the node removal is that the quality of the solution might decrease since the search region for new paths shrinks each time paths are fixed due to the network reduction. Note that this only holds true if the optimal solution contains empty trips, which is rare in practice.

\section{Results}\label{section:results}
We now study the performance of our heuristics using a real-world data set. First, we describe the data set and explain how we derive instances of varying size from it. To evaluate the performances of the heuristics, we compare the heuristics by using fixed parameter values on some smaller instances. Thereafter, we tune the parameters of the best-performing heuristic on the smaller instances. Lastly, we solve the larger instances.

The heuristics are implemented in Python 3.7 including the CPLEX solver version 12.10.0. We use a MacBook Pro with a 1.4 GHz Quad-Core Intel Core i5 processor and 8 GB RAM.

\subsection{Case}\label{section:data}
We now describe the data set used for the evaluation of our methodology. The E-VSP requires input on the timetabled trips, the depots, the charging locations, and the bus types. We use data on the timetable of the bus concession of Gooi en Vechtstreek provided by Lynxx to test our proposed heuristics. The data set is used to create several instances with a varying number of trips. 

The bus concession of Gooi en Vechtstreek covers the public bus transit in Bussum, Hilversum, Huizen, Naarden, and Weesp \citep{OVinNederland2021Concessie2011-2021}. These cities are all located in the eastern part of the province of Noord-Holland. We consider the timetable of December the 3rd, 2019. In total, 816 trips are scheduled that day. The information regarding the trips comes from a General Transit Feed Specification (GTFS) data set. GTFS includes information about the public transit trips and corresponding geographic information. 
 
The trips in this data set can be divided in several clusters. The \textit{city lines} contain only stops that are located in Hilversum. The trips from the city lines are on average shorter than the trips from the \textit{Rnet line}, which contains trips between Amsterdam and Hilversum. Trips corresponding to the \textit{regional lines} depart from various places. The \textit{rush hour lines} and \textit{school lines} contain 27 and 3 trips, respectively. These lines are only used at specific times. 

When viewing the number of trips operated simultaneously during the day, two peak moments occur. After the start of the service day, we see an increasing number of trips until a peak in the morning around 08:00. Afterward, the number of trips decreases, stabilizes, and increases again after approximately 15:00. The evening peak lasts until approximately 18:30. Thereafter, the number of trips decreases. The maximum number of trips that should be serviced simultaneously is 53. This means that at least 53 buses are required for the bus schedule of the entire concession.

The data set includes two depots, three charging stations, and two bus types. Both depots have a large capacity. Two of the charging stations are located close to the depots, and have a capacity of five and two buses, respectively. The third charging station has a capacity of two buses.

\begin{table}
    \caption{Bus characteristics of two bus types that are used in the case to provide a bus schedule for the concession of Gooi en Vechtstreek} \label{tab:bus_characteristics}
    \centering
    \scalebox{0.75}{
    \begin{tabular}{|c|c|c|}
    \hline
    \textbf{Characteristics}    & \textbf{Type 1} & \textbf{Type 2}  \\ \hline
    Battery Capacity (kWh) & 155 & 210    \\
    Consumption Rate (kWh/km) & 1.3   & 1.4  \\ 
    Idling Consumption Rate (kWh/s) & 0.00167 & 0.00167\\
    Charging Rate (kWh/s) & 0.0639 & 0.0889 \\
    Investment Costs (\EUR{}) & 50,000 & 52,500 \\ 
    Operational Costs (\EUR{}/km) & 1.0 & 1.05 \\ \hline
    \end{tabular}}
\end{table}

Finally, we consider two different bus types. Information about the bus types regarding their battery capacities, consumption and charging rates, and costs is given in Table~\ref{tab:bus_characteristics}. We translate the battery capacity to a value of $100\%$ for the SoC level. Then, we convert the information that is given in kWh in Table~\ref{tab:bus_characteristics} to SoC values in proportion to the battery capacity.
In this case, we assume both bus types can be located at both depots, can operate all trips and can charge at every charging station.

\subsection{Parameter Settings and Instances}\label{section:instances&param}
We now first describe the parameter settings for our heuristics. Then, we use the data set described in Section~\ref{section:data} to create instances with varying numbers of trips. 

\subsubsection{Parameter Settings}
We initially fix the discretization of the possible SoC values and the length of the time blocks. These fixed parameter values are used to compare the heuristics. In Section~\ref{section:sensitivity_analysis}, we study the impact of coarser and finer discretizations. The discretization influences the number of nodes and arcs in the pricing networks. We use the same discretization for all networks. 

To compare the heuristics,  we use 27 possible SoC values, ranging from 22\% to 100\% in steps of 3\%. We start from 22\% to incorporate a minimum value for the SoC.
Furthermore, we choose 5 minutes as the length of time blocks. We assume that the charging rate is independent of the SoC value at the beginning of the charging action, as it is often assumed in literature \citep[see e.g.][]{VanKootenNiekerk2017SchedulingVehicles}. Note that this could be easily adjusted because of the way we construct the network. 

When creating the network for each instance, we set the maximum deadhead time to 1 hour. The maximum time allowed for idling is set to 8 hours between trips and to 3 hours if a bus goes to or comes from a charging station. We determine the deadhead distance and times using The Open Source Routing Machine (OSRM). In the data set that we consider, all deadhead times are less than 1 hour.
Moreover, we use energy costs of €0.1361 per kWh and crew costs of €0.67 per minute, similar to \citet{VanAken2020StrategischeOptimizer}. We use €10 as starting costs of a charging activity. Other cost parameters are given in Table~\ref{tab:bus_characteristics}

\subsubsection{Instances}\label{section:subsectionInstances}
Several instances are used for the computational results. First, we create relatively smaller instances by picking random subsets of the set of all trips. Secondly, we use the clusters as described in Section~\ref{section:data} to create different instances. Characteristics of the instances are shown in Table~\ref{tab:instances}. Here, Trips is the number of trips in the instance, TH is the length of the time horizon, and Nodes and Arcs represent the total number of nodes and arcs in the four networks, respectively.  
The time horizon, used for the charging nodes, depends on the length of the service for the considered trips. For each instance, we let the time horizon begin at the start of the hour of the first trip that departs and end at the end of the hour of the last trip that finishes.
Instance~A represents 50 trips randomly chosen from all morning trips, containing all trips that finish before noon. Instance~B represents 100 trips randomly chosen from all trips. Instances 1 to 5 are combinations of the clusters. 

\begin{table}[h]
    \caption{Characteristics of the instances} \label{tab:instances}
    \centering
    \scalebox{0.75}{
    \begin{tabular}{|c|c|D{.}{.}{3.0}|D{.}{.}{2.0}|D{.}{.}{7.0}|D{.}{.}{10.0}|}
    \hline
    \textbf{Instance}    & \textbf{Trips used} & \multicolumn{1}{c|}{\textbf{Trips (\#)}} & \multicolumn{1}{c|}{\textbf{TH (h)}} & \multicolumn{1}{c|}{\textbf{Nodes (\#)}} &  \multicolumn{1}{c|}{\textbf{Arcs (\#)}  }\\ \hline
    A & Random Morning Trips & 50 & 6\rlap{$^1$} & 13,691\rlap{$^1$} & 224,873\rlap{$^1$}\\ 
    B & Random Trips & 100 &  20\rlap{$^1$} &68,345\rlap{$^1$}  & 1,229,069\rlap{$^1$} \\ 
    1 & City line & 119 & 19 & 70,507 & 2,084,615 \\ 
    2 & Rnet line & 185 & 21 & 74,671 & 1,339,960 \\ 
    3 & City \& Rnet lines & 304 & 21 & 86,846 & 3,800,870 \\ 
    4 & $\begin{array}{c}
         \textnormal{Rush hour- \& School-}  \\
          \textnormal{\& Regional Lines}
    \end{array}$ & 512 & 21 &103,134  & 9,170,356 \\ 
    5 & All bus lines & 816 & 21 & 125,344 & 15,589,878 \\ \hline
    \multicolumn{6}{l}{$^1$ Average outcome of ten instances}
    \end{tabular}
		}
\end{table}

In general, the number of nodes and arcs increases if the number of trips increases. However, we see that instance 1 contains more arcs in the final networks than instance 2. This can be explained by the fact that instance 1 includes the shorter City line trips, in contrast to the Rnet line trips that have a longer average duration in instance 2. This results in fewer arcs for instance 2 than for instance 1, since fewer trip pairs are compatible.   

\subsection{Comparison Heuristics}\label{section:comparison_heuristics}
We now compare the performance of the heuristics we propose in Section~\ref{section:meth-heuristics}. We first state the parameter values we use. Second, we show the computational results of the comparison.

\subsubsection{Heuristic Parameters}
For the truncated price-and-branch and the truncated column generation heuristics, we stop the column generation process if there is no improvement of at least 0.01\% ($Z_{\textnormal{min}}=0.01$) within 30 iterations ($I=30$). 
We use a larger number of iterations than \citet{Pepin2009AProblem} to prevent early termination due to degeneracy. 
Additionally, since we use a relative decrease and the objective value is large at the beginning of the algorithm, a larger value of $I$ prevents the algorithm from terminating too soon. The variables with a value higher than 0.70 are fixed in the truncated column generation heuristic ($\theta = 0.70$), similar to \citet{Pepin2009AProblem}. For the (truncated) price-and-branch heuristic, we set a time limit of 1 hour for CPLEX for solving the BP.

\subsubsection{Results}\label{section:allHeuristics}
In this section, we test the proposed heuristics on the smaller instances A, B, 1, 2, and 3. For instances A and B, the results are averaged over 10 randomly generated instances. The required time and final solution values of the heuristics are presented in Table~\ref{tab:comparisonheur}. For all heuristics, Time is the total time needed for the heuristic to obtain a final integer solution, including the time needed to solve the BP for the (truncated) price-and-branch heuristic. Sol is the objective value of the final integer solution of the E-VSP, and G is a bound on the optimality gap. 
The bounds on the optimality gap of all heuristics are computed using a lower bound that is obtained using an optimistic rounding scheme, as explained in Section~\ref{section:lowerBounds}.
The values of these lower bounds and more detailed results can be found in Appendix~\ref{appendix:detailedresults}.

\begin{table}[h]
\caption{A comparison of the required time and final solution of the three proposed heuristics} \label{tab:comparisonheur}
\centering
\scalebox{0.75}{
\begin{tabular}{c|D{.}{.}{6.0}D{.}{.}{9.0}D{.}{.}{6.0}|D{.}{.}{6.0}D{.}{.}{9.0}D{.}{.}{6.0}|D{.}{.}{6.0}D{.}{.}{9.0}D{.}{.}{4.0}}
\cline{1-10}
\textbf{} & \multicolumn{3}{c|}{\textbf{Price-and-Branch}} & \multicolumn{3}{c|}{\textbf{Truncated Price-and-Branch}} & \multicolumn{3}{c}{\textbf{Truncated Column Generation}} \\
\textbf{} & \multicolumn{3}{c|}{} & \multicolumn{3}{c|}{($Z_{\textnormal{min}}=0.01, I=30$)} & \multicolumn{3}{c}{($Z_{\textnormal{min}}=0.01, I=30, \theta = 0.70$)} \\ \hline
Instance  & \multicolumn{1}{c}{Time (s)}&  \multicolumn{1}{c}{Sol}  &\multicolumn{1}{c|}{ G (\%)} &\multicolumn{1}{c}{ Time (s)}  & \multicolumn{1}{c}{Sol}  &\multicolumn{1}{c|}{ G (\%)}  & \multicolumn{1}{c}{Time (s)}& \multicolumn{1}{c}{Sol}&\multicolumn{1}{c}{ G (\%)}\\ \hline
A  & $35$\rlap{$^1$} & $744,828$\rlap{$^1$} & $0.869$\rlap{$^1$}& $26$\rlap{$^1$} & $744,896$\rlap{$^1$} & $0.878$\rlap{$^1$} & $37$\rlap{$^1$} & $744,838$\rlap{$^1$} & $0.870$\rlap{$^1$} \\
B  & $1,436$\rlap{$^1$} & $569,503$\rlap{$^1$} & $0.149$\rlap{$^1$} & $451$\rlap{$^1$} & $602,477$\rlap{$^1$} & $6.165$\rlap{$^1$} & $962$\rlap{$^1$} & $569,614$\rlap{$^1$} & $0.167$\rlap{$^1$} \\
1 & $9,504$ &  $304,320$ & $19.901$ & 792 &  304,320 & $19.901$  & 1,977 & 254,525 & $0.282$ \\
2  & $6,163$\rlap{$^2$} &  $1,153,237$\rlap{$^2$} & $6.027$\rlap{$^2$} & $4,394$\rlap{$^2$}  & $1,319,157$\rlap{$^2$} & $21.281$\rlap{$^2$} & 2,274 &  1,102,030 & $1.319$ \\
 3  & $22,204$\rlap{$^2$} &  $1,549,954$\rlap{$^2$} & $20.278$\rlap{$^2$}   & $6,603$\rlap{$^2$} &  $1,773,327$\rlap{$^2$} & $37.612$\rlap{$^2$}  & 10,134 & 1,299,192 & $0.818$ \\\hline
\multicolumn{6}{l}{$^{1}$ Average outcome of ten instances} \\
\multicolumn{6}{l}{$^{2}$ BP not solved to optimality due to reached time limit}
\end{tabular}}
\end{table}

Table~\ref{tab:comparisonheur} shows that all three heuristics obtain a high-quality solution for the A-instances. Both the price-and-branch and the truncated column generation heuristic provide a good solution to the B-instances as well. For the larger instances 1, 2, and 3, we see that the truncated column generation outperforms the other two heuristics based on solution quality. 

The price-and-branch heuristic requires more than twice as much computation time as the truncated column generation for instances 2 and 3. Note that for instances 2 and~3, the time limit for solving the BP is reached for both the price-and-branch and the truncated price-and-branch heuristic. 
The price-and-branch heuristic also requires most time for Instance 1 and the B-instances. Despite its longer computation time, the price-and-branch heuristic provides a poor optimality gap bound of about 20\% for instance 1. The truncated price-and-branch heuristic is the fastest for all instances except instance 2. However, the truncated price-and-branch heuristic results in the solution with the highest costs for all instances and has a poor optimality gap bound ($>19\%$) for the larger instances.

The required time can be partly explained by the number of iterations needed in the column generation process of all heuristics. We present in Table~\ref{tab:iterations_heur} the required number of iterations and the average time for solving the pricing problem and RMP per iteration. Here, It is the number of required iterations during the column generation process, PP is the average time of solving the pricing problem per iteration, and RMP is the average time of solving the RMP per iteration.

\begin{table}[h]
\caption{The number of iterations and the average time per iteration for solving the pricing problem and the RMP of the three proposed heuristics} \label{tab:iterations_heur}
\centering
\scalebox{0.75}{
\begin{tabular}{c|D{.}{.}{5.0}D{.}{.}{4.1}D{.}{.}{4.0}|D{.}{.}{4.1}D{.}{.}{7.4}D{.}{.}{5.1}|D{.}{.}{5.0}D{.}{.}{7.4}D{.}{.}{5.1}}
\cline{1-10}
\textbf{} & \multicolumn{3}{c|}{\textbf{Price-and-Branch}} & \multicolumn{3}{c|}{\textbf{Truncated Price-and-Branch}} & \multicolumn{3}{c}{\textbf{Truncated Column Generation}} \\
\textbf{} & \multicolumn{3}{c|}{} & \multicolumn{3}{c|}{($Z_{\textnormal{min}}=0.01, I=30$)} & \multicolumn{3}{c}{($Z_{\textnormal{min}}=0.01, I=30, \theta = 0.70$)} \\ \hline
Instance  & \multicolumn{1}{c}{It (\#)}& \multicolumn{1}{c}{PP (s)}  & \multicolumn{1}{c|}{RMP (s)} & \multicolumn{1}{c}{It (\#)}& \multicolumn{1}{c}{PP (s)}  & \multicolumn{1}{c|}{RMP (s)}  & \multicolumn{1}{c}{It (\#)}& \multicolumn{1}{c}{PP (s)}  & \multicolumn{1}{c}{RMP (s)} \\ \hline
A  & 154\rlap{$^1$} & $0.22$\rlap{$^1$} & $0.004$\rlap{$^1$} &  114\rlap{$^1$} &	$0.22$\rlap{$^1$} &$0.004$\rlap{$^1$} & 169\rlap{$^1$} & $0.22$\rlap{$^1$} & $0.004$\rlap{$^1$}  \\
B  &$1,076$\rlap{$^1$} & $1.24$\rlap{$^1$} & $0.011$\rlap{$^1$} &$258$\rlap{$^1$} & $1.25$\rlap{$^1$} & $0.010$\rlap{$^1$} &$767$\rlap{$^1$} & $1.26$\rlap{$^1$} & $0.010$\rlap{$^1$}  \\
1  &  2,947  & $2.12$ & $0.016$ & 369 & $2.13$ & $0.010$  & 936 & $2.11$ & $0.011$  \\
 2 & 1,835 & $1.37$ & $0.015$ & 580 & $1.33$ & $0.013$   &  1,673 & $1.35$ & $0.014$ \\
 3  &  4,601 & $3.97$ & $0.056$    & 754 & $3.95$ & $0.016$  & 2,657 & $3.81$ & $0.021$ \\\hline
\multicolumn{6}{l}{$^1$ Average outcome of ten instances}
\end{tabular}}
\end{table}

As expected, the price-and-branch heuristic requires more iterations than the truncated price-and-branch heuristic due to the early stopping criterion of the truncated price-and-branch heuristic. Similarly, the truncated column generation heuristic requires more iterations than the truncated price-and-branch heuristic, since the truncated column generation heuristic continues the column generation process after the root node has been explored, in contrast to the truncated price-and-branch. We observe that the pricing problem dominates the computation time per iteration, and that the average time to solve one pricing problem is highly similar for the three heuristics. Thus, the computation time is determined mainly by the number of iterations.

Based on the above-described results and the trade-off between computation time and solution quality, we conclude that the truncated column generation heuristic is the best-performing heuristic. Hence, we continue with the truncated column generation heuristic in the remainder of this section.

\subsection{Parameter tuning}\label{section:sensitivity_analysis}
In this section, we perform a parameter tuning to study the influence of the chosen parameter values for the truncated column generation heuristic, using the three smallest instances (A, B, and 1). In contrast to the previous section, we now consider a single instance A and B. First, we solve the instances using different values for the parameters of the truncated column generation heuristic to test their influence on the solution quality. Afterward, we solve the instances using various discretizations for the underlying network to test how this impacts the final solution.

 \subsubsection{Heuristic Parameters}\label{section:sensitivity-heuristic}
We test how the quality of the solution changes if other parameters are used for the truncated column generation heuristic. We fix the discretization and use possible SoC values in steps of 3\% between 22\% and 100\% and time blocks of 5 minutes, similar to Section~\ref{section:comparison_heuristics}. We test different values for the minimum relative improvement ($Z_{\textnormal{min}}$), the number of iterations ($I$), and the threshold for fixing variables ($\theta $) used in the truncated column generation heuristic. Note that the minimum allowed value for $\theta$ is 0.5. Smaller values could lead to violations of the capacity constraints. 

The results for Instance~1 are presented in Table~\ref{tab:sensitivityparameters}. Here, Time is the total computation time of the heuristic, It is the number of iterations, and PP and RMP are the average times per iteration to solve the pricing problem and the RMP, respectively. Sol, G, and B are the objective value, a bound on the optimality gap, and the number of buses of the final integer solution, respectively. For instance A and B, we obtained highly similar results for all settings. For this reason, we have reported the results for these instances in the appendix.

\begin{table}[h]
 \caption{Results of the truncated column generation heuristic for differing parameters considering the minimum required relative improvement, the number of iterations, and the threshold as used in the fixing step} \label{tab:sensitivityparameters}
 \centering
\scalebox{0.75}{
\begin{tabular}{D{.}{.}{4.0}D{.}{.}{2.0}D{.}{.}{2.0}|D{.}{.}{4.0}D{.}{.}{5.0}D{.}{.}{5.0}D{.}{.}{4.0}D{.}{.}{9.2}D{.}{.}{6.2}D{.}{.}{2.3}D{.}{.}{2.0}}
\hline
  \multicolumn{3}{c|}{\textbf{}} & \multicolumn{8}{l}{\textbf{Instance 1} ($\textnormal{LB}= 253,809.6$)} \\ \cline{4-11}
 \multicolumn{1}{c}{$Z_\textnormal{min} (\%)$}  & \multicolumn{1}{c}{$I$}  &
 \multicolumn{1}{c|}{$\theta$}  &  \multicolumn{1}{c}{Time (s)} &  \multicolumn{1}{c}{It (\#)} &\multicolumn{1}{c}{PP (s)}  & \multicolumn{1}{c}{RMP (s)} & \multicolumn{1}{c}{Sol} & \multicolumn{1}{c}{G (\%)}  & \multicolumn{1}{c}{ B (\#) }& \\ \hline
 $0.010$ & 30  & $0.7$ & 1,977 & 936 & $2.11$ & $0.011$ & 254,525 & $0.282$ & 5  \\ 
 $0.010$ & 15 & $0.7$ & 44 & 21 & $2.20$ & $0.008$ & 304,320 & $19.901$ & 6  \\
 $0.010$ & 50  & $0.7$ & 3,309 & 1,514 & $2.17$ & $0.012$ & 254,327 & $0.204$  & 5  \\ 
  $0.010$ & 90  & $0.7$ & 4,234 & 1,942 & $2.16$ & $0.013$ & 254,250 & $0.173$  & 5  \\
 $0.005$ & 30  & $0.7$ & 3,359 & 1,554 & $2.15$ & $0.013$ & 254,287 & $0.188$ & 5  \\ 
$0.050$ & 30  & $0.7$ &  1,912 & 897 & $2.13$ & $0.011$ & 254,661 & $0.335$ & 5 \\
$0.500$ & 30  & $0.7$ & 1,854 & 890 & $2.08$ & $0.010$ & 304,737 & $20.065$ & 6 \\ 
$0.010$ & 30  & $0.5$ & 1,978 & 936 & $2.11$ & $0.011$ & 254,525 & $0.282$ & 5   \\
$0.010$ & 30  & $0.9$  & 1,960 & 936 & $2.09$ & $0.011$ & 254,525 & $0.282$ & 5 \\\hline
\end{tabular}}
\end{table}

Table~\ref{tab:sensitivityparameters} demonstrates that using a value for $I$ of 15 results in a worse solution than using larger values for $I$ (30, 50, or 90). This can partly be explained by the occurrence of degeneracy, which results in the objective value not improving for several iterations. This causes the column generation process to be stopped too early if $I$ is small. In general, increasing $I$ results in a better solution. However, it also increases the number of required iterations and, as a consequence, the required computation time. 

Decreasing the value of $Z_\textnormal{min}$ to 0.005 while keeping $I$ and $\theta$ constant results in slightly better solutions compared to using $Z_\textnormal{min}=0.01$. However, more iterations are executed, which again increases the computation time. Increasing $Z_\textnormal{min}$ from 0.01 to 0.05 or 0.5 results in a lower computation time, but also in an increased solution value.  

The usage of different values for $\theta$ does not always result in a difference in the final solutions.
This holds if none of the path variables are above the threshold of 0.7 when the early stopping criterion is met. In this case, the path variable with the maximum value is fixed. This is confirmed by the results for varying values of $\theta$: all values we have tested result in the same solution.

\subsubsection{Discretization} \label{section:results-discretization}
As introduced in Section~\ref{section:literature_review}, finding a balance between better solutions versus computational tractability plays a crucial role when determining the best discretization \citep{Boland2019TheDesign}. We now study the impact of the discretization on the performance of the heuristic. In these tests, the heuristic parameters are kept constant with $Z_{\textnormal{min}}=0.01$, $I=30$, and $ \theta = 0.70$, similarly as in Section~\ref{section:comparison_heuristics}. Table~\ref{tab:sensitivitydiscretization} presents the results of the experiments. We test different values for the length of the time blocks (TB) and the step size of the SoC values (Steps). Range represents the interval of the possible SoC values.

\begin{table}[h]
 \caption{Results of the truncated column generation heuristic for different discretization levels.} \label{tab:sensitivitydiscretization}
 \centering
\scalebox{0.75}{
\begin{tabular}{D{.}{.}{2.0}D{.}{.}{2.0}c|D{.}{.}{5.0}D{.}{.}{5.0}D{.}{.}{4.0}D{.}{.}{5.0}D{.}{.}{9.0}D{.}{.}{2.3}D{.}{.}{2.0}}
\hline
\multicolumn{3}{c|}{\textbf{}} & \multicolumn{7}{l}{\textbf{Instance A} ($\textnormal{LB}= 653,319.3$)} \\ \cline{4-10}
 \multicolumn{1}{c}{TB (min)}  & \multicolumn{1}{c}{Steps (\%)}  &
 \multicolumn{1}{c|}{Range}  &  \multicolumn{1}{c}{Time (s)} &  \multicolumn{1}{c}{It (\#)} &\multicolumn{1}{c}{PP (s)}  & \multicolumn{1}{c}{RMP (s)} & \multicolumn{1}{c}{Sol}  & \multicolumn{1}{c}{G (\%)} & \multicolumn{1}{c}{B (\#) } \\ \hline
 5 &3  & 22-100 & 45 & 181 & $0.25$ & $0.004$ & 655,957 & 0.404 & 13   \\
  2 &3  & 22-100 & 125 & 205 & $0.61$ & $0.006$ & 655,986 & 0.408 & 13  \\
  10 &3  & 22-100 & 25 & 210 & $0.12$  & $0.003$ & 705,972 & 8.059 & 14  \\
30 &3  & 22-100 &8 & 206 & $0.04$  & $0.002$ & 658,618 & 0.811 & 13  \\
1 & 1 & 22-100 &931 & 213 & $4.44$ & $0.012$ & 655,901 & 0.395 &  13 \\
  5 & 1  & 22-100 & 168 & 236 & $0.71$  & $0.004$ & 655,928 & 0.399 & 13  \\
  5 & 6  & 22-100 & 25 & 191  & $0.12$  & $0.004$ & 655,991 & 0.409 & 13  \\
5 &10  & 20-100 & 14 & 191  & $0.07$ & $0.003$ & 703,452 & 7.674 & 14   \\
 5 &20  & 20-100 & 3 & 120 & $0.02$ & $0.003$ & 1,236,037 & 89.193 & 24  \\\hline
\multicolumn{3}{c|}{\textbf{}} & \multicolumn{7}{l}{\textbf{Instance B} ($\textnormal{LB}= 559,350.2$)} \\ \cline{4-10}
 \multicolumn{1}{c}{TB (min)}  & \multicolumn{1}{c}{Steps (\%)}  &
 \multicolumn{1}{c|}{Range}  &  \multicolumn{1}{c}{Time (s)} &  \multicolumn{1}{c}{It (\#)} &\multicolumn{1}{c}{PP (s)}  & \multicolumn{1}{c}{RMP (s)} & \multicolumn{1}{c}{Sol}  & \multicolumn{1}{c}{G (\%)} & \multicolumn{1}{c}{B (\#) } \\ \hline
  5 &3  & 22-100 &  850 & 700 & $1.21$ & $0.010$ & 559,863 & 0.092 & 11  \\ 
  2 &3  & 22-100 & 2,301 & 801 & $2.87$ & $0.017$ & 560,079 & 0.130 & 11 \\
  10 &3  & 22-100 & 470 & 762 & $0.61$  & $0.007$ & 560,371 & 0.183 & 11    \\
30 &3  & 22-100 & 165 & 796 & $0.20$  & $0.004$ & 564,990 & 1.008 & 11  \\
  5 &1   & 22-100 & 3,320 & 766 & $4.37$  & $0.011$ & 559,776 & 0.076 & 11  \\
  5 & 6  & 22-100 & 366 & 707  & $0.51$  & $0.009$ & 560,107 & 0.135 & 11  \\
5 &10  & 20-100 & 220 & 801  & $0.26$ & $0.009$ & 562,810 & 0.619 & 11   \\
 5 &20  & 20-100 &24 & 283 & $0.08$ & $0.008$ & 2,623,753 & 369.072 & 51  \\\hline
 \multicolumn{3}{c|}{\textbf{}} & \multicolumn{7}{l}{\textbf{Instance 1} ($\textnormal{LB}= 253,809.6$)} \\ \cline{4-10}
 \multicolumn{1}{c}{TB (min)}  & \multicolumn{1}{c}{Steps (\%)}  &
 \multicolumn{1}{c|}{Range}  &  \multicolumn{1}{c}{Time (s)} &  \multicolumn{1}{c}{It (\#)} &\multicolumn{1}{c}{PP (s)}  & \multicolumn{1}{c}{RMP (s)} & \multicolumn{1}{c}{Sol}  & \multicolumn{1}{c}{G (\%)} & \multicolumn{1}{c}{B (\#) } \\ \hline
  5 &3  & 22-100 & 1,977 & 936 & $2.11$ & $0.011$ & 254,525  & 0.282 & 5   \\
  2 &3  & 22-100 &3,941 & 1,175 & $3.32$ & $0.018$ & 254,555 & 0.294 & 5   \\
  10 &3  & 22-100 & 45 & 37 & $1.25$ & $0.009$ & 304,725  & 20.060 & 6   \\
30 &3  & 22-100 & 616 & 1,143 & $0.53$ & $0.007$ & 304,639  & 20.027 & 6   \\
  5 &1   & 22-100 & 8,444 & 1,078 & $7.85$ & $0.012$ & 254,459 & 0.256 & 5   \\
  5 & 6  & 22-100 &  34 & 37 & $0.93$ & $0.008$ & 304,670 & 20.039 & 6   \\
5 &10  & 20-100 & 18 & 37 & $0.50$ & $0.008$ & 304,991 & 20.165 & 6   \\
 5 & 20  & 20-100 & 76 & 458 & $0.16$ & $0.009$ & 2,003,821 & 689.498 & 40  \\\hline
    \end{tabular}}
\end{table}

First, note that time blocks of 5 minutes and a SoC step size of 3\% results in small optimality gaps of 0.40\%, 0.09\% and 0.28\% for instance A, B and 1, respectively. Since these gaps are obtained using a lower bound that does not depend on the used discretization, only relatively small improvements are theoretically possible if one uses a finer discretization. Indeed, we find that choosing a different length of the time blocks than 5 minutes while keeping the possible SoC step size constant at 3\%, does rarely result in a better solution. Using time blocks of 2 minutes instead of 5 minutes even worsens the solution quality, despite requiring more iterations and increasing the computation time. This could be explained by the amount of charging in the considered time block and the conservative rounding. In 5 minutes, the buses charge approximately 12.5\%. In 2 minutes, the increase in SoC is roughly 5\%. Using a step size of 3\%, the amount that is discarded due to rounding between two consecutive charging actions is approximately 0.5\% versus 2\%, for time blocks of 5 versus 2 minutes, respectively. Hence, using time blocks of 2 minutes might result in a stronger underestimation of the SoC. This demonstrates the interplay between the discretization of time and that of the SoC values. In particular, it is not guaranteed that using a finer time discretization results in a better solution. The usage of a finer time discretization does result in a longer computation time, though.

An increased computation time can also be observed for time blocks of 1 minute and steps of 1\% for the possible SoC values for instance A. The cost decrease is negligible, but the computation time is increased by a factor 20. To avoid computation times that are prohibitively long, the experiments with time blocks of 1 minute and steps of 1\% for the possible SoC values are not executed for instances B and 1. 

Table~\ref{tab:sensitivitydiscretization} also shows the influence of choosing different step sizes than 3\% for the SoC values for a time block length of 5 minutes. For instances A, B, and 1, choosing a step size of 1\% results in bus schedules with slightly better solution values, for example due to fewer scheduled charging actions. 
However, this also approximately triples the average time needed to solve the pricing problem per iteration. This increases the required computation time significantly. 

For instances A and B, a slightly worse solution is obtained by using a step size of 6\% instead of 3\%. 
Using a step size of 6\% roughly halves the time needed to solve the pricing problem per iteration. This causes a significant reduction in the required computation time. We note that a step size of 6\% results in the same amount of SoC discarded between two consecutive charging actions as a step size of 3\%. Using a step size of 10\% causes a solution that requires one more bus for instance A. 

For all instances, using a step size of 20\% results in a solution that requires significantly more buses. This can be explained by the fact that both buses can increase their SoC by approximately 12.5\% during a charging action of 5 minutes. 
Thus, in this setting, charging actions do actually not result in an increase in the SoC using this network structure. Hence, choosing a step size for the possible SoC values larger than the amount of SoC increase after charging one time block results in solutions with high costs and requiring more buses. 

To conclude, using time blocks of 5 minutes in combination with a step size of 3\% gives the best solutions compared to using different lengths of the time blocks. In general, using a step size of 1\% results in slightly better solutions than using 3\%, but also increases the required computation time. On the other hand, enlarging the step size to 6\% shortens the average time to solve the pricing problem per iteration. However, it also can slightly worsen the solution value.

\subsection{Larger Instances} \label{section:larger_instances}
In this section, we solve the two larger instances. First, we use instance 4 and the results of Section~\ref{section:sensitivity_analysis} to compare different discretized values for this larger instance using the truncated column generation heuristic. Additionally, we also test the truncated column generation with node removal heuristic. Afterward, we solve the whole concession using the insights obtained for instance~4.

\subsubsection{Instance 4}\label{section:instance4}

In this section, we solve instance 4 using various techniques.
We do not deviate from the used parameter values in Section~\ref{section:comparison_heuristics} for the truncated column generation heuristic. Hence, we use $Z_{\textnormal{min}}=0.01$,  $I=30$, and $ \theta = 0.70$. Using these parameter settings, we obtained good solutions for instances A, B, 1, 2, and 3 as shown in Section~\ref{section:comparison_heuristics}. All achieved optimality gap bounds are less than 1.5\%. 

For the discretized values, we solve instance 4 using time blocks of 5 minutes and a step size of 3\% between 22\% and 100\% for the SoC values. Based on the results in Section~\ref{section:results-discretization}, the time blocks of 5 minutes give the best solutions. 
Besides, we prefer using a step size of 3\% because of the lower computation times compared to a step size of 1\%, despite the slightly worse solutions. 
 
To obtain a solution in less time, we also solve instance~4 using a step size of 6\% instead of 3\% for the possible SoC values. Doing so leads to smaller networks, and allows the pricing problem to be solved faster.
Additionally, we test the node removal heuristic to further reduce the average time needed to solve the pricing problem. This extension removes specific nodes and arcs during the column generation process each time paths are fixed. 

\begin{table}[]
\caption{Results of the truncated column generation (with and without node removal) heuristic for instance~4 using  $Z_{\textnormal{min}}=0.01$, $I=30$, $\theta = 0.70$, TB = 5 min.\ and Range = 22-100. The lower bound for this instance equals 1,528,574.} \label{tab:instance4tcg}
 \centering
\scalebox{0.75}{

\begin{tabular}{cc|rrrrrrr}
\hline
\multicolumn{1}{l}{Steps (\%)} & \multicolumn{1}{l|}{Node Removal} & \multicolumn{1}{l}{Time (h)} & \multicolumn{1}{l}{It (\#)} & \multicolumn{1}{l}{PP (s)} & \multicolumn{1}{l}{RMP (s)} & \multicolumn{1}{l}{Sol} & \multicolumn{1}{l}{G (\%)} & \multicolumn{1}{l}{B (\#)} \\ \hline
3                              & No                                & 26.00                        & 9,348                       & 9.87                       & 0.13                        & 1,598,917               & 4.60                       & 31                         \\
6                              & No                                & 10.66                        & 8,555                       & 4.35                       & 0.12                        & 1,622,100               & 6.12                       & 31                         \\
3                              & Yes                               & 13.83                        & 8,100                       & 5.96                       & 0.13                        & 1,543,099               & 0.95                       & 30                         \\
6                              & Yes                               & 7.20                         & 9,472                       & 2.56                       & 0.13                        & 1,573,957               & 2.97                       & 30                         \\ \hline
\end{tabular}           
}
\end{table}

The results are presented in Table~\ref{tab:instance4tcg}. We use the same abbreviations as in  Table~\ref{tab:sensitivityparameters}. We observe that the truncated column generation heuristic using a step size of 3\% requires 26 hours of computation time to solve instance 4.  A large part of this computation time is needed to solve the pricing problem, which takes on average 9.87 seconds per iteration. This can be explained by the large number of nodes and arcs in the corresponding network (see Table~\ref{tab:instances}).

Using a step size of 6\%, the time needed to solve the pricing problem is significantly reduced to 4.35 seconds on average.
This is more than two times faster than solving the pricing problem using a step size of 3\%. However, we also see an increase in the solution value of €23,183, which is an increase of 1.45\%. 
The total time needed for the column generation process is decreased to less than 11 hours. 

Using the truncated column generation with node removal heuristic, we see a further decrease in the average time needed to solve the pricing problem. Now, on average 5.96 and 2.56 seconds are needed to solve the pricing problem per iteration, for steps of 3\% and 6\%, respectively. Thus, the node removal heuristic significantly reduces the computation time. In particular, 
 using steps of 6\%, the total computation time is reduced to 7 hours. 

We also observe a decrease in the solution value if nodes are removed from the network. We believe this is due to the heuristic nature of the solution method. 

\subsubsection{Instance 5: The Entire Concession}

Based on the results of Section~\ref{section:instance4}, we use the truncated column generation with node removal heuristic to solve the entire concession (instance 5) with 816 trips. For the heuristic parameters we use $Z_{\textnormal{min}}=0.01$, $I=30$, and $\theta = 0.70$. For the underlying network, we use time blocks with a length of 5 minutes and SoC values between 22\% and 100\% with a step size of 6\%. This results in a network containing in total 63,035 nodes and 7,487,553 arcs. Note that this is a significant reduction compared to the number of nodes and arcs reported in Table~\ref{tab:instances}.

\begin{table}[h]
 \caption{Results of the truncated column generation with node removal heuristic for instance~5 using  $Z_{\textnormal{min}}=0.01$, $I=30$, and $\theta = 0.70$ } \label{tab:instance5_extended}
 \centering
\scalebox{0.75}{
\begin{tabular}{ccc|cccccccc}
\hline
\multicolumn{1}{c}{TB (min)}  & \multicolumn{1}{c}{Steps (\%)}  &
 \multicolumn{1}{c|}{Range}   & \multicolumn{1}{c}{Time (h)} &  \multicolumn{1}{c}{It (\#)}& \multicolumn{1}{c}{PP (s)}  & \multicolumn{1}{c}{RMP (s)} & \multicolumn{1}{c}{Sol} & \multicolumn{1}{c}{LB} & \multicolumn{1}{c}{G (\%)}  & B (\#) \\ \hline
5 & 6 & 22-100  & 28.00 & 20,266 & 4.30 & 0.59 & 2,794,568 & 2,702,417 & 3.41 & 53 \\ \hline
\end{tabular}}
\end{table}
We present the results of solving the entire concession in Table~\ref{tab:instance5_extended}. 
The truncated column generation with node removal heuristic requires 28 hours to solve instance 5. The optimality gap is below 3.5\%. The optimal solution value of the RMP equals 2,755,589. Thus, roughly half of the optimality gap can be explained by the discretization, and the other half by the heuristic to obtain integer solutions.
The feasible bus schedule requires 53 buses. This equals the minimum number of buses required based on the number of trips that are serviced simultaneously as explained in Section~\ref{section:data}. In total, 824 trips are scheduled. This means that 8 trips are driven `empty'. 

Below, we discuss a few characteristics of the bus duties that comprise the final solution. 
Each bus charges around 2 hours on average during its duty. The average times of deadheading and idling per bus duty are approximately 79 and 166 minutes, respectively. Due to the usage of discretized SoC values in our methodology, the SoC is rounded down in the resulting bus duties before a charging action or before servicing a trip. Per bus duty, the amount of SoC that is discarded as a result adds up to 66 percentage points on average. If we track the SoC of the bus duties without rounding down,
we see that the average minimum real SoC value is 35.7\%, with an overall minimum of 27.5\%. This is higher than the minimum allowed SoC value of 22\%. In contrast, if the SoC is rounded down, all bus duties reach the minimum allowed SoC value of 22\% at some point in time. The higher actual minimum SoC values show that the resulting bus duties are feasible, but possibly also suboptimal, since they do not exploit the entire range of allowed SoC values in reality.
 
\begin{figure}[ht]
    \centering
    \includegraphics[scale=0.35]{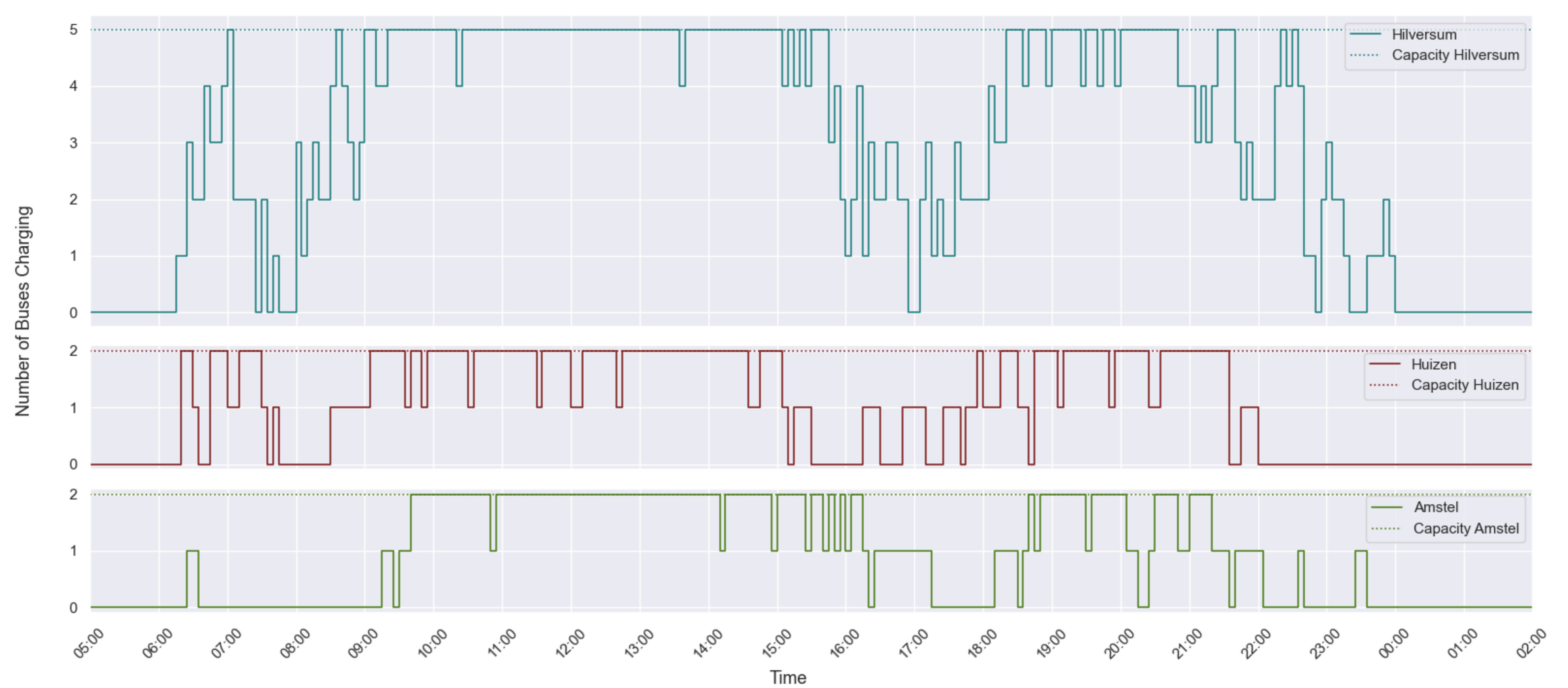}
    \caption{The occupation of the charging stations in the resulting bus schedule}
    \label{fig:capacities}
\end{figure}

An important aspect of a feasible bus schedule is that the capacities of the charging stations are never exceeded. We illustrate the usage of the charging stations during the day in Figure~\ref{fig:capacities}. 
Because we incorporated the charging capacities in our methodology, the capacities of the charging stations are never exceeded. This can be seen in Figure~\ref{fig:capacities}. Additionally, we see that the charging stations are often fully occupied,
especially after the morning peak and evening peak hours. In contrast, during these peak moments, the number of buses that are charging is low.

\section{Conclusion and Future Research}\label{section:conclusion}
In this paper, we studied the Electric Vehicle Scheduling Problem considering both the capacity of charging stations and partial charging. To solve the problem, we developed a network structure in which every path represents a feasible duty, and presented two heuristics based on column generation. Valid lower bounds were obtained by solving an auxiliary problem on a slightly altered network. 

Computational results based on a bus concession in the Netherlands showed that the truncated column generation outperforms price-and-branch, achieving an optimality gap below 1.5\% on the smaller instances. When applied in combination with the node removal heuristic to speed up the pricing problem, truncated column generation solved the entire concession with 816 trips to an optimality gap of 3.4\%. 

There are numerous promising directions for future research. It is likely that, in addition to the node removal heuristic, the pricing problem can be sped up further by developing other dedicated acceleration strategies. It would also be interesting to jointly optimize the vehicle schedule with the driver schedule, as the possibility to coordinate charging actions and meal breaks introduces interesting dynamics into the problem. Finally, our method to compute true lower bounds irrespective of the discretization could serve as the starting point for a full-fledged dynamic discretization discovery algorithm, potentially finding better solutions. 

\vspace{10pt}
\noindent 
{ \footnotesize
\textbf{Acknowledgements} We would like to thank the company Lynxx for facilitating and supervising the project that inspired this paper. Marelot H. de Vos did most of the work on this project while she was at Erasmus School of Economics, Erasmus University Rotterdam.
}

\addcontentsline{toc}{section}{References}
\bibliography{References}
\bibliographystyle{apacite}

\clearpage
\appendix
\section{Detailed Description of Arcs}\label{appendix:arcs}
In this section, we provide a detailed description of the arcs that are present in the primal and in the dual network. 

Consider a given combination $k\in \mathcal{K}$ of a vehicle type and a depot. In the primal network, we apply a conservative rounding scheme. Here, we define the function $F^\textnormal{soc}(\cdot)$ that returns the largest SoC value $s \in \mathcal{S}^{k}$ smaller than or equal to the SoC input. In the dual network, we apply an optimistic rounding scheme. There, the function $F^\textnormal{soc}(\cdot)$ is redefined and returns the smallest SoC value $s \in \mathcal{S}^{k}$ that is larger than or equal to the SoC input.

The compatibility of two nodes  $n \in \mathcal{N}^{k}$ and $v \in \mathcal{N}^{k}$ depends on several characteristics. Recall that $\tau_{(n,v)}^\textnormal{soc}$ denotes the SoC required for deadheading and idling between nodes $n$ and $v$. Similarly, we define $\phi_{(n,v)}^\textnormal{time}$ and $\chi_{(n,v)}^\textnormal{time}$ as the required deadheading and idling time between nodes $n$ and $v$, respectively. Moreover, let $\phi_{(n,v)}^\textnormal{max}$ and $\chi_{(n,v)}^\textnormal{max}$ denote the maximum allowed deadhead and idle times between nodes $n$ and $v$. These quantities depend on the combination $n$ and $v$ to be able to impose different restrictions per location. The parameters $b_i$ and $e_i$ denote the begin and end time of trip $i \in \mathcal{T}$. 

Below, we explain which arcs are created in network $k$ while discussing the outgoing arcs from the source node, the trip nodes, and the charging nodes separately. Accordingly, we define the sets $A_{d^{\sigma}_{k}}$, $\mathcal{A}_{\textnormal{trip}_k}$, and $A_{\textnormal{charge}_k}$. 
The complete set of arcs of network $k$ is then given by $$\mathcal{A}^{k}=A_{d^{\sigma}_{k}} \cup \mathcal{A}_{\textnormal{trip}_k} \cup A_{\textnormal{charge}_k}.$$

\paragraph{Source Node} First, we consider the source node $d^{\sigma}_{k}$. By our assumption that vehicles leave the depot fully charged, we do not include arcs that go from the source node to a charging node. When creating outgoing arcs to trip nodes $v \in \mathcal{N}_{k}^\textnormal{trip}$, we consider the required SoC usage for the deadheading from the depot to the start location of the corresponding trip. An arc is created if and only if the maximum deadhead time is not exceeded. Given that the vehicle can depart from the depot such that it arrives precisely in time to operate the trip, idling need not be considered. We obtain the set $$A_{d^{\sigma}_{k}}=\{(n,v)|n=d^{\sigma}_{k}, v=(i,s) \in \mathcal{N}_{k}^\textnormal{trip},s=F^\textnormal{soc}(s^\textnormal{full}-\tau_{(n,v)}^\textnormal{soc}), \phi_{(n,v)}^{\textnormal{time}} \leq \phi_{(n,v)}^{\textnormal{max}}\}.$$

\paragraph{Trip Node} Next, we explain the outgoing arcs from trip nodes $n\in \mathcal{N}_{k}^\textnormal{trip}$ by considering outgoing arcs to the sink node, other trip nodes, and charging nodes separately, defining the sets $\mathcal{A}_{\textnormal{trip}_k}^\textnormal{sink}$, $\mathcal{A}_{\textnormal{trip}_k}^\textnormal{trip}$, and $\mathcal{A}_{\textnormal{trip}_k}^\textnormal{charge}$, respectively. 

From trip node $n=(i,s)$, an outgoing arc to the sink node is created if the corresponding SoC value $s$ suffices to return to the depot after performing the trip. Here, we consider the SoC required for the trip, for the deadheading back to the depot, and take into account that the SoC value should always stay above $s_k^\textnormal{min}$. This results in the set $$\mathcal{A}_{\textnormal{trip}_k}^\textnormal{sink}= \{(n,v)|n=(i,s) \in \mathcal{N}_{k}^{\textnormal{trip}}, v=d^{\tau}_{k},s \ge f_i +\tau_{(n,v)}^\textnormal{soc} + s_k^\textnormal{min},\phi_{(n,v)}^{\textnormal{time}} \leq \phi_{(n,v)}^{\textnormal{max}}\}.$$ 

Outgoing arcs from trip node $n=(i,s)$ to other trip nodes $v=(j,s') \in \mathcal{N}_{k}^\textnormal{trip}$ are only possible if the two considered trips are compatible while not exceeding the maximum deadhead and idle times. The idle time is defined as $\chi^\textnormal{time}_{(n,v)}=b_j - e_i - \phi^\textnormal{time}_{(n,v)}$. We denote $i \shortrightarrow j$ for the requirements $\phi^\textnormal{time}_{(n,v)} \le \phi_{(n,v)}^\textnormal{max}$ and $0\leq \chi^\textnormal{time}_{(n,v)}\le \chi_{(n,v)}^\textnormal{max}$. Additionally, from the SoC $s$ of trip node $n$, the SoC required for the trip corresponding to node $n$ and for the deadheading and idling between the nodes is subtracted to determine the SoC $s'$  of the successor node $v \in \mathcal{N}_{k}^\textnormal{trip}$. This results in the set
\begin{equation*}
\mathcal{A}_{\textnormal{trip}_k}^\textnormal{trip}= \{(n,v)|n=(i,s) \in \mathcal{N}_{k}^{\textnormal{trip}}, v=(j,s')\in \mathcal{N}_{k}^\textnormal{trip},
i \shortrightarrow j, s' = F^\textnormal{soc}(s - f_{i} - \tau^\textnormal{soc}_{(n,v)}) \}.
\end{equation*}

 A charging node $v=(r,b,s') \in \mathcal{N}_{k}^\textnormal{charge}$ can directly succeed the trip node $n$ if $\phi^\textnormal{time}_{(n,v)} \le \phi_{(n,v)}^\textnormal{max}$ and $0 \leq \chi^\textnormal{time}_{(n,v)} \le \chi_{(n,v)}^\textnormal{max}$, where we now compute the idle time as $\chi^\textnormal{time}_{(n,v)} = t_b - e_i - \phi^\textnormal{time}_{(n,v)}$. Thus, the begin time $t_b$ of time block $b$ should be later than the end time of trip $i$ plus the deadheading time from the end location of trip $i$ to the location of the charging station $r$. Additionally, the maximum allowed deadhead and idle time cannot be exceeded. For brevity, we denote these requirements by $i\shortrightarrow (r,b)$. For the SoC value $s'$ of the charging node $v$, the required SoC for the trip of node $n$ as well as the SoC needed for the deadheading and idling from the end location of the trip to the charging station is subtracted from the SoC value $s$ belonging to node $n$. This results in the arcs from the set
\begin{equation*}
\mathcal{A}_{\textnormal{trip}_k}^\textnormal{charge}= \{(n,v)|n=(i,s) \in \mathcal{N}_{k}^{\textnormal{trip}}, v=(r,b,s') \in \mathcal{N}_{k}^\textnormal{charge}, i \shortrightarrow (r,b), s' = F^\textnormal{soc}(s - f_{i} - \tau^\textnormal{soc}_{(n,v)}) \}.
\end{equation*}

In the dual network, we want to allow a vehicle to start charging immediately when arriving at the charging station. For that reason, if $0\leq \chi^\textnormal{time}_{(n,v)}< l$, we add the amount of SoC that can be charged during idling. In particular, the last requirement is replaced by $s' = F^\textnormal{soc}(s + \hat{s} - f_{i} - \tau^\textnormal{soc}_{(n,v)})$, where $\hat{s}$ is the amount of SoC that could be charged while the vehicle is idling, and where $F^\textnormal{soc}$ now uses an optimistic rounding scheme. In practice, this corresponds to the situation where the vehicle starts charging during the time block before time block $b$, immediately after arriving. The vehicle is considered not to occupy the charging station during this time block. 

The union of the sets $\mathcal{A}_{\textnormal{trip}_k}^\textnormal{sink}$, $\mathcal{A}_{\textnormal{trip}_k}^\textnormal{trip}$, and $\mathcal{A}_{\textnormal{trip}_k}^\textnormal{charge}$ provides the set $\mathcal{A}_{\textnormal{trip}_k}$ that includes all outgoing arcs from the trip nodes $n \in \mathcal{N}_{k}^\textnormal{trip}$.

\paragraph{Charging Node} Lastly, we discuss the outgoing arcs from each charging node $n \in \mathcal{N}_{k}^\textnormal{charge}$ by again considering outgoing arcs to the sink node, trip nodes, and other charging nodes separately. Recall that $s_n^\textnormal{+}$ represents the SoC increase during the charging action corresponding to charging node $n$. The amount of increase depends on the length of the time block of charging node $n$ and the charging rate. 

Since each vehicle can charge at its depot, an outgoing arc from charging node~$n$ to the sink node is only created if the corresponding SoC value without the charging of node $n$ is too low for deadheading back to the depot, while it is high enough after the charging. Hence, returning to the depot from a charging station is only possible if it is necessary to charge before returning to the depot. This results in the set $$\mathcal{A}_{\textnormal{charge}_k}^\textnormal{sink}= \{(n,v)|n=(r,t,s)\in \mathcal{N}_{k}^\textnormal{charge}, v=d^{\tau}_{k},s < \tau_{(n,v)}^\textnormal{soc} + s_k^\textnormal{min} \le s +  s_n^\textnormal{+}, \phi_{(n,v)}^{\textnormal{time}} \leq \phi_{(n,v)}^{\textnormal{max}}\}.$$ 

A trip node $v=(i,s')\in  \mathcal{N}_{k}^\textnormal{trip}$ can directly succeed charging node $n=(r,b,s)$ if $\phi^\textnormal{time}_{(n,v)} \le \phi_{(n,v)}^\textnormal{max}$ and $0 \leq \chi^\textnormal{time}_{(n,v)} \le \chi_{(n,v)}^\textnormal{max}$, where we compute the idle time as $\chi^\textnormal{time}_{(n,v)} = b_i - t_b - l - \phi^\textnormal{time}_{(n,v)}$. Thus, the begin time of the trip is later than the end time of the time block corresponding to node $n$ plus the required deadhead time, and the maximum allowed deadhead and idle times are not exceeded. We denote these requirements by $(r,b)\shortrightarrow i$. The SoC value $s'$ of the trip node $v$ should be equal to the highest possible SoC value smaller than the SoC value $s$ of the charging node $n$ plus the increase due to the charging and minus the required SoC for the deadheading and idling. This results in the set of outgoing arcs
\begin{equation*}
\mathcal{A}_{\textnormal{charge}_k}^\textnormal{trip}= \{(n,v)|n=(r,b,s)\!\in\! \mathcal{N}_{k}^\textnormal{charge}, v=(i,s') \!\in\! \mathcal{N}_{k}^\textnormal{trip}, (r,b) \shortrightarrow i , s' = F^\textnormal{soc}(s +s_n^+ - \tau^\textnormal{soc}_{(n,v)})\}.
\end{equation*}

In the dual network, the arc $(n,v)$ is present if the vehicle can execute the trip by departing during time block $b$, instead of at the end of time block $b$. In particular, $0 \leq \chi^\textnormal{time}_{(n,v)}$ is then replaced by $-l > \chi^\textnormal{time}_{(n,v)}$. If $\chi^\textnormal{time}_{(n,v)}<0$, the vehicle charges only during the time it is actually present at the charging station. Thus, $s_n^+$ is then replaced by $\hat{s}$, where $\hat{s}$ is the amount of SoC that can be charged during the time $l+\chi^\textnormal{time}_{(n,v)}$.

We assume that a charging activity cannot be paused and cannot continue at a different location. Accordingly, an outgoing arc from charging node $n$ to another charging node $v~\in~\mathcal{N}_{k}^\textnormal{charge}$ is created only if both nodes have the same charging station and the time block of charging node $v$ begins immediately when the time block of charging node $n$ ends. Due to this, a charging activity can take multiple time blocks, but it is not allowed for a vehicle to switch between different charging stations or for the charging activity to contain a break. Additionally, an outgoing arc to charging node $v$ is created only if the value of the SoC increases. 
This results in the set
$$\mathcal{A}_{\textnormal{charge}_k}^\textnormal{charge}= \{(n,v)|n=(r,b,s)\in\mathcal{N}_{k}^\textnormal{charge},v=(r',b',s') \in \mathcal{N}_{k}^\textnormal{charge},n\rightarrow v,s'=F^\textnormal{soc}(s+s^+_n)>s\},$$
where $n\rightarrow v$  denotes the requirements $r = r'$, $t_{b'} =t_b+l$.

The union of the sets $\mathcal{A}_{\textnormal{charge}_k}^\textnormal{sink}$, $\mathcal{A}_{\textnormal{charge}_k}^\textnormal{trip}$, and $\mathcal{A}_{\textnormal{charge}_k}^\textnormal{charge}$ provides the set $A_{\textnormal{charge}_k}$ including all outgoing arcs from the charging nodes.

\section{Detailed Results Comparison Heuristics}\label{appendix:detailedresults}

Table~\ref{tab:details_comparison_heur_app} presents the detailed results of the runs that are discussed in Section~\ref{section:allHeuristics}. We first report the computation times. For the (truncated) price-and-branch heuristics, Time CG is the time needed for the column generation process and Time BP is the time needed to solve the binary program. For the truncated column generation, Time is the total time needed to solve the instances. Moreover, 
It is the number of iterations of the column generation process, PP is the average time for solving the pricing problem per iteration, and RMP is the average time for solving the RMP per iteration. LB is a lower bound on the optimal solution value. Sol is the objective value of the E-VSP corresponding to the final integer solution of the heuristic, G is the optimality gap bound, and B is the number of buses used in the final bus schedule. 

\begin{table}[ht]
 \caption{Detailed results of all three proposed heuristics} \label{tab:details_comparison_heur_app}
 \centering
\scalebox{0.75}{
\begin{tabular}{c|rrrrrrrrr}
\hline
\multicolumn{1}{c|}{\textbf{}} & \multicolumn{9}{l}{\textbf{Price-and-Branch}} \\ \cline{2-10}
 \multicolumn{1}{c|}{Instance}  &  \multicolumn{1}{r}{Time CG (s)} & \multicolumn{1}{r}{Time BP(s)} & \multicolumn{1}{r}{It (\#)} & \multicolumn{1}{r}{PP (s)} & \multicolumn{1}{r}{RMP (s)}  & \multicolumn{1}{r}{LB} & \multicolumn{1}{r}{Sol} & \multicolumn{1}{r}{G (\%)} & \multicolumn{1}{r}{ B (\#) }\\ \hline
 A & 35.20\rlap{$^1$} & 0.02\rlap{$^1$}& 154.40\rlap{$^1$} &	0.22\rlap{$^1$} & 0.004\rlap{$^1$} & 738,971.06\rlap{$^1$} &	744,828.45\rlap{$^1$} &	0.869\rlap{$^1$} & 14.80\rlap{$^1$} \\ 
B & $1,356.78$\rlap{$^1$} & $78.79$\rlap{$^1$} & $1,076.10$\rlap{$^1$} & $1.24$\rlap{$^1$} & $0.011$\rlap{$^1$} & $568,726.21$\rlap{$^1$} &	$569,502.93$\rlap{$^1$} &	$0.149$\rlap{$^1$} & $11.20$\rlap{$^1$} \\
1  &  6,346.01  & 3,158.03 & 2,947.00  & 2.12 & 0.016 & 253,809.57 & 304,320.39 & 19.901 & 6.00 \\ 
2 & 2,555.74 & $3,607.29$\rlap{$^2$} & 1,835.00 & 1.37 & 0.015 & 1,087,685.65 & 1,153,236.86 & 6.027 & 22.00 \\
3 & 18,604.00 & $3,600.10$\rlap{$^2$} & 4,601.00 & 3.97 & 0.056 & 1,288,645.26 & 1,549,954.37 & 20.278 & 30.00 \\ \hline
\multicolumn{1}{c|}{\textbf{}} & \multicolumn{9}{l}{\textbf{Truncated Price-and-Branch ($Z_{\textnormal{min}}=0.01, I=30$)}} \\ \cline{2-10}
\multicolumn{1}{c|}{Instance}  & \multicolumn{1}{r}{Time CG (s)} & \multicolumn{1}{r}{Time BP(s)} & \multicolumn{1}{r}{It (\#)} & \multicolumn{1}{r}{PP (s)} & \multicolumn{1}{r}{RMP (s)} & \multicolumn{1}{r}{LB} & \multicolumn{1}{r}{Sol}& \multicolumn{1}{r}{G (\%)}     &\multicolumn{1}{r}{ B (\#) }  \\ \hline
A& 25.72\rlap{$^1$}	& 0.04\rlap{$^1$} &	113.50\rlap{$^1$} &	0.22\rlap{$^1$} & 0.004\rlap{$^1$} & 738,971.06\rlap{$^1$} & 744,895.65\rlap{$^1$} &	0.878\rlap{$^1$} & 14.80\rlap{$^1$}  \\ 
B & $326.29$\rlap{$^1$} &	$125.19$\rlap{$^1$} & $258.40$\rlap{$^1$} & $1.25$\rlap{$^1$}& $0.010$\rlap{$^1$} & $568,726.21$\rlap{$^1$} & $602,477.38$\rlap{$^1$} &	$6.165$\rlap{$^1$} & $11.70$\rlap{$^1$} \\ 
 1 & 791.20 & 0.98 & 369.00 & 2.13 & 0.010 & 253,809.57 & 304,320.39 & 19.901 & 6.00 \\ 
 2 & 782.47 & $3,611.40$\rlap{$^2$} & 580.00 & 1.33 & 0.013 & 1,087,685.65 & 1,319,156.97 & 21.281 & 25.00  \\
3 & 2,998.59 & $3,604.73$\rlap{$^2$} & 754.00 & 3.95 & 0.016 & 1,288,645.26 & 1,773,326.53 & 37.612 & 34.00 \\ \hline
\multicolumn{1}{c|}{\textbf{}} & \multicolumn{9}{l}{\textbf{Truncated Column Generation ($Z_{\textnormal{min}}=0.01, I=30, \theta = 0.70$)}} \\ \cline{2-10}
\multicolumn{1}{c|}{Instance}  & \multicolumn{1}{r}{Time (s)} & \multicolumn{1}{r}{} &  \multicolumn{1}{r}{It (\#)}& \multicolumn{1}{r}{PP (s)}  & \multicolumn{1}{r}{RMP (s)} & \multicolumn{1}{r}{LB} & \multicolumn{1}{r}{Sol} & \multicolumn{1}{r}{G (\%)}   &\multicolumn{1}{r}{ B (\#) }  \\ \hline 
A & 37.33\rlap{$^1$} && 168.90\rlap{$^1$} & 0.22\rlap{$^1$} &	0.004\rlap{$^1$}& 738,971.06\rlap{$^1$}& 744,837.94\rlap{$^1$} &	0.870\rlap{$^1$} & 14.80\rlap{$^1$} \\ 
B &$962.00$\rlap{$^1$} &&	$766.90$\rlap{$^1$} & $1.26$\rlap{$^1$} &	$0.010$\rlap{$^1$} & $568,726.21$\rlap{$^1$} & $569,613.50$\rlap{$^1$} & $0.167$\rlap{$^1$} & $11.20$\rlap{$^1$} \\ 
1 & 1,976.95 && 936.00 & 2.11 & 0.011 & 253,809.57 & 254,524.52 & 0.282 & 5.00 \\ 
2 & 2,273.95 && 1,673.00 & 1.35 & 0.014 & 1,087,685.65 & 1,102,030.14 & 1.319 & 21.00  \\
3 & 10,134.22 && 2,657.00 & 3.81 & 0.021 & 1,288,645.26 & 1,299,191.73 & 0.818 & 25.00 \\ \hline
\multicolumn{6}{l}{$^1$ Average outcome of 10 instances} \\
\multicolumn{6}{l}{$^{2}$ BP not solved to optimality due to reached time limit}
\end{tabular}}
\end{table}

Table~\ref{tab:sensitivityparametersAll} presents the  results that are obtained for various heuristic parameter settings on Instance A and Instance B, as discussed in Section~\ref{section:sensitivity-heuristic}  Here, Time is the total computation time of the heuristic, It is the number of iterations, and PP and RMP are the average times per iteration to solve the pricing problem and the RMP, respectively. Sol, G, and B are the objective value, a bound on the optimality gap, and the number of buses of the final integer solution, respectively.

\begin{table}[ht]
 \caption{Results of the truncated column generation heuristic for differing parameters considering the minimum required relative improvement, the number of iterations, and the threshold as used in the fixing step} \label{tab:sensitivityparametersAll}
 \centering
\scalebox{0.75}{
\begin{tabular}{D{.}{.}{4.0}D{.}{.}{2.0}D{.}{.}{2.0}|D{.}{.}{4.0}D{.}{.}{5.0}D{.}{.}{5.0}D{.}{.}{4.0}D{.}{.}{9.2}D{.}{.}{6.2}D{.}{.}{2.3}D{.}{.}{2.0}}
\hline
\multicolumn{3}{c|}{\textbf{}} & \multicolumn{8}{l}{\textbf{Instance A} ($\textnormal{LB}= 653,319.3$)} \\ \cline{4-11}
 \multicolumn{1}{c}{$Z_\textnormal{min} (\%)$}  & \multicolumn{1}{c}{$I$}  &
 \multicolumn{1}{c|}{$\theta$}  &  \multicolumn{1}{c}{Time (s)} &  \multicolumn{1}{c}{It (\#)} &\multicolumn{1}{c}{PP (s)}  & \multicolumn{1}{c}{RMP (s)} & \multicolumn{1}{c}{Sol} & \multicolumn{1}{c}{G (\%)}  & \multicolumn{1}{c}{ B (\#) }& \\ \hline
 $0.010$ & 30  & $0.7$ & 46 & 181 & $0.25$ & $0.004$ & 655,957 & $0.404$ & 13  \\ 
 $0.010$ & 15 & $0.7$ & 53 & 215 & $0.25$ & $0.004$ & 656,012 & $0.412$ & 13  \\
 $0.010$ & 50  & $0.7$ & 44 & 171 & $0.25$ & $0.004$ & 655,941 &  $0.401$ & 13  \\
 $0.010$ & 90  & $0.7$ & 47 & 186 & $0.25$ & $0.004$ & 655,946 & $0.402$  & 13 \\
 $0.005$ & 30  & $0.7$ & 48 & 191 & $0.25$ & $0.004$ & 655,947 & $0.402$ & 13 \\
 $0.050$ & 30  & $0.7$ & 55 & 221 & $0.25$ & $0.004$ & 655,950 & $0.402$ & 13 \\
$0.500$ & 30  & $0.7$ & 48 & 189 & $0.26$ & $0.004$ & 655,961 & $0.404$ & 13 \\ 
 $0.010$ & 30  & $0.5$ & 42 & 162 & $0.25$ & $0.004$ & 655,965 & $0.405$ & 13 \\
 $0.010$ & 30  & $0.9$ & 43 & 172 & $0.25$ & $0.004$ & 655,941 & $0.401$ & 13 \\\hline
 \multicolumn{3}{c|}{\textbf{}} & \multicolumn{8}{l}{\textbf{Instance B} ($\textnormal{LB}= 559,350.2$)} \\ \cline{4-11}
 \multicolumn{1}{c}{$Z_\textnormal{min} (\%)$}  & \multicolumn{1}{c}{$I$}  &
 \multicolumn{1}{c|}{$\theta$}  &  \multicolumn{1}{c}{Time (s)} &  \multicolumn{1}{c}{It (\#)} &\multicolumn{1}{c}{PP (s)}  & \multicolumn{1}{c}{RMP (s)} & \multicolumn{1}{c}{Sol} & \multicolumn{1}{c}{G (\%)}  & \multicolumn{1}{c}{ B (\#) }& \\ \hline
 $0.010$ & 30  & $0.7$ & 850 & 700 & $1.21$ & $0.010$ & 559,863 & $0.091$ & 11  \\ 
 $0.010$ & 15 & $0.7$ & 612 & 513 & $1.20$ & $0.009$ & 560,638 & $0.230$ & 11  \\
 $0.010$ & 50  & $0.7$ & 1,199 & 975 & $1.22$ & $0.010$ & 559,819 & $0.083$  & 11  \\ 
 $0.010$ & 90  & $0.7$ & 1,570 & 1,277 & $1.22$ & $0.011$ & 559,635 & $0.051$  & 11  \\ 
 $0.005$ & 30  & $0.7$ & 1,108 & 884 & $1.25$ & $0.010$ & 559,775 & $0.076$ & 11  \\
 $0.050$ & 30  & $0.7$ &  495 & 416 & $1.20$ & $0.009$ & 560,379 & $0.184$ & 11 \\
 $0.500$ & 30  & $0.7$ & 424 & 361 & $1.19$ & $0.009$ & 560,350 & $0.178$ & 11 \\ 
 $0.010$ & 30  & $0.5$ & 1,006 & 807 & $1.24$ & $0.010$ & 559,852 & $0.089$ & 11   \\
 $0.010$ & 30  & $0.9$ & 893 & 700 & $1.38$ & $0.010$ & 559,863 & $0.091$ & 11 \\\hline
\end{tabular}}
\end{table}

\end{document}